\documentclass[11pt]{article}
\author{Hao Chen\\
Department of Mathematics\\
Hangzhou Dianzi University\\
Hangzhou 310018, Zhejiang Province, China\\
E-mail: haochen@hdu.edu.cn}
\title{ A ternary construction of lattices}
\date{July, 2014}

\begin{document}

\maketitle

\begin{abstract}
In this paper we propose a general ternary construction of lattices from three rows and ternary codes. Most laminated lattices and Kappa lattices in ${\bf R}^n$, $n\leq 24$ can be recovered from our tenary construction naturally. This ternary construction of lattices can  be used to generate many new "sub-optimal"  lattices of low dimensions. Based on this ternary construction of lattices new extremal even lattices of dimensions $32,40$ and $48$  are also constructed.

\end{abstract}

\section{Introduction}

How one can arrange most densely in space an infinite number of equal spheres is a classical mathematical problem and a part of
Hilbert 18th problem (\cite{CS1,Co}). It is deeply rooted in information theory and physics (\cite{Torquato,Kallus,MT}). The root
lattices in Euclid spaces of dimensions 1, 2, 3, 4, 5, 6, 7 and 8
and the Leech lattice of dimension $24$ (\cite{CS1}) had been proved  to be the unique densest lattice sphere packings
in these dimensions (see \cite{CS1,Mar,CoKu}). For a historic survey we refer to  \cite{CS1} Chapter 1 and \cite{Mar}.
From Voronoi's theory (\cite{Mar,Sch}), there are algorithms to determine the densest lattice sphere packings for all dimensions.
However the computational task for dimensions $n > 9$ is generally infeasible. Laminated lattices in ${\bf R}^n, n \leq24$,
were known from the work of A. Korkine, G.Zolotareff in 1877, T. W. Chaundy in 1946 and J . Leech in 1967 (\cite{CS1}, page 158-159, \cite{Mar}).
In 1982 J. H. Conway and N. J. A. Sloane calculated all densities of laminated lattices
up to dimension 48 (\cite{CS2}). The $n$ dimensional laminated lattice, $1 \leq n \leq 24$ and $n \neq 11, 12, 13$,
is the only known densest lattice sphere packing in ${\bf R}^n$ (\cite{CS1}, page 15). The known densest lattice in
${\bf R}^n$, $n=11,12,13$ is the Kappa lattice ${\bf K}_{11}$, the Coxeter-Todd lattice ${\bf K}_{12}$ and the Kappa
lattice ${\bf K}_{13}$(\cite{CS1}, page 15). \\

The lattice ${\bf
L}^{*}=\{{\bf y} \in {\bf R}^m: <{\bf y},{\bf x}> \in {\bf Z}, {\bf x} \in {\bf L}\}$ is
called the dual lattice of the lattice ${\bf L}$. A lattice is
called integral if the inner products between lattice vectors are
integers and an integral lattice is called even if the Euclid norms of all lattice vectors are even numbers. An integral lattice satisfying ${\bf L}={\bf L}^{*}$ is called unimodular lattice. The minimum norm $\mu({\bf L})$ for an unimoular even lattice satisfies  $\mu({\bf L}) \leq 2[\frac{n}{24}]+2$ and an unimodular even lattice satisfying the equality is called an extremal even unimodular lattice (\cite{CS1}). The most interesting cases are the jumping dimensions $n \equiv 0$ $mod$ $24$. In dimension $24$ J. H. Conway proved that the Leech lattice is the unique extremal even unimodular lattice {(\cite{Conway}). In dimension $48$, there are four constructed extremal even unimodular lattices (\cite{CS1,Nebe}. It has long been an open problem  if there is an extremal even extremal even unimodualr lattice in dimension $72$.  The recent two breakthroughs are the first consctructed $72$ dimensional extremal even unimodular lattice (\cite{Griess,Nebe1}) and the fourth $48$ dimensional extremal even unimodular lattice (\cite{Nebe2}). We refer to \cite{MOS,Bachoc,Gaborit,Harada,Harada1,HKMV,Gross} for the general bounds and the general construction of extremal even or odd unimodular or modular lattices. It is a classical result of C. L. Mallows, A. M. Odlyzko and N. J. A. Sloane that there is no $n$ dimensional extremal even unimodular lattice when $n$ is sufficiently large (\cite{MOS}). \\

Though it is well-known that for the dimension $32$, there are millions of extremal even unimodular lattices (\cite{King}), there are only $15$ known  constructed $32$ dimensional extremal even unimodular lattices found by Koch and Venkov (\cite{KochVen,KochNebe,Nelist}).  We refer to  \cite{Calderbank,Ozeki,Nebe}, \cite{CS1} page 221 for the extremal even unimodular lattices of dimension $40$ and  \cite{BachocNebe,StehleWatkins,Watkins} for the four constructed $80$ dimensional extremal even unimodular lattices.\\

In this paper we propose a construction of lattices from ternary
codes.  The ternary construction recovers the known densest
 lattices including the Leech lattice ${\bf \Lambda}_{24}$,  the Coxeter-Todd lattice ${\bf K}_{12}$, the Barnes-Wall lattice ${\bf \Lambda}_{16}$,  and one of the laminated lattice ${\bf \Lambda}_{26}$ of dimension $26$ naturally.  Many other laminated and Kappa lattices in dimensions $n \leq 24$ can also be recovered from our ternary construction naturally.
We generate some new "sub-optimal" low dimensional lattices from
this ternary construction. Based on this ternary construction of lattices new extemal even lattices of dimensions $32,40,48$ are also constructed.\\

For a packing of equal non-overlapping spheres
in ${\bf R^n}$ with centers ${\bf x_1},{\bf x_2}, ...$, ${\bf
x_m},....$, the packing radius $\rho$ is $\frac{1}{2}min_{i \neq
j}\{||{\bf x_i}-{\bf x_j}||\}$. The density $\Delta$ is  $lim_{t
\rightarrow 0} \frac{Vol\{{\bf x} \in {\bf R^n}: ||{\bf x}|| <t,\
\exists {\bf x_i}, ||{\bf x}-{\bf x_i}||<\rho\}}{Vol\{{\bf x} \in
{\bf R^n}: ||{\bf x}||<t\}}$. The center density $\delta$ is
$\frac{\Delta}{V_n}$ where $V_n$ is the volume of the ball of
radius $1$ in ${\bf R^n}$. Let ${\bf b}_1,...,{\bf b}_m$ be $m$ linearly independent vectors in
the Euclid space ${\bf R}^n$ of dimension $n$. The discrete point
sets ${\bf L}=\{x_1{\bf b}_1+\cdots+x_n{\bf b}_m: x_1,...,x_m \in
{\bf Z}\}$ is a dimension $m$ lattice in ${\bf R}^n$.  The volume of the
lattice is $Vol({\bf L})=(det(<{\bf b}_i, {\bf
b}_j>))^{\frac{1}{2}}$. Let $\lambda({\bf L})$ be the length of
the shortest non-zero vectors in the lattice and the minimum norm of
the lattice is just $\mu({\bf L})=(\lambda({\bf L}))^2$. The set of lattice vectors
with the length $\lambda({\bf L})$ is denoted by $min({\bf L})$ and the cardinality of
$min({\bf L})$ is the kissing number $k_1({\bf L})$ of the lattice ${\bf L}$. Let $\mu'({\bf L})$ be the minimum norm of lattices vectors in the set ${\bf L} \setminus \{{\bf 0}\} \setminus min({\bf L})$. The number of the second layer lattice vectors $k_2({\bf L})$ is the number of lattice vectors with the norm $\mu'({\bf L})$. When the
centers of the spheres are these lattice vectors in ${\bf L}$ we
have a lattice sphere packing for which $\rho=\frac{1}{2}\lambda({\bf L})$ and the center density
$\delta({\bf L})=\frac{\rho^n}{Vol({\bf L})}$.\\

Let $q$ be a prime power and ${\bf F}_q$ be the finite field with
$q$ elements. A linear (non-linear ) error-correcting code ${\bf C}
\subset {\bf F}_q^n$ is a $k$ dimensional subspace(or a subset of
$M$ vectors). For a codeword ${\bf x} \in {\bf C}$, the support $Supp({\bf x})$ of ${\bf x}$ is the set of
positions of non-zero coordinates. The Hamming weight $wt({\bf x})$ is
the number of elements in the support of ${\bf x}$. The minimum Hamming
weight(or distance) of the linear (or non-linear) code $C$ is
defined as $d({\bf C})=min_{{\bf x}\neq {\bf y}, {\bf x},{\bf y} \in {\bf C}}\{wt({\bf x-y})\}$.
We refer to $[n,k,d]_q$ (or $(n,M,d)_q$) code as linear (or
non-linear)code with length $n$,  distance $d$ and dimension $k$ (or
$M$ codewords). For a binary code ${\bf C} \subset {\bf F}_2^n$ the
construction A (\cite{CS1}) leads to a lattice in ${\bf R}^n$.
The lattice ${\bf L(C)}$ is defined as the set of integral vectors ${\bf
x}=(x_1,...,x_n) \in {\bf Z}^n$ satisfying $x_i\equiv c_i$ $mod$ $2$
for some codeword ${\bf c}=(c_1,...,c_n) \in {\bf C}$. This is a lattice with the center density
$\delta=\frac{min\{\sqrt{ d({\bf C})},2\}^n}{2^{2n-k({\bf C})}}$.
This construction A gives some best known densest lattice packings
in low dimensions(see \cite{CS1}). For a $(n,M,d)$ non-linear binary code, the
same construction gives the non-lattice packing with center density
$\frac{M \cdot min\{\sqrt{d({\bf C})},2\}^n}{2^{2n}}$. Some of
them are the known best sphere packings (see \cite{CS1}) or have presently known highest kissing numbers. The known densest packing (non-lattice) in dimension 10 with center density $\frac{5}{128}$ is from non-linear binary $(10, 40, 4)$ code found by M. R. Best in 1980 (\cite{CS1}). The non-lattice packing ${\bf P}_{11a}$ in ${\bf R}^{11}$ with the center density $\frac{9}{256}$ and the kissing number $566$ is constructed from a non-linear binary $(11, 72, 4)$ code (\cite{CS1}, page 139). \\

\section{The main construction}

{\bf Lemma 2.1.} {\em Let ${\bf L}$ be a dimension $r$ lattice in ${\bf R}^n$ with the volume $vol({\bf L})$ and the center density $\delta({\bf L})$. Let ${\bf E}$ be the dimension $2r$ lattice in ${\bf R}^{3n}$ defined by ${\bf E}=\{({\bf x},{\bf y},{\bf z}): {\bf x} \in {\bf L}, {\bf y} \in {\bf L}, {\bf z} \in {\bf L}, {\bf x}+{\bf y}+{\bf z}=0 \}$. Then the volume of the ${\bf E}$ is $3^{\frac{r}{2}} \cdot vol({\bf L})^2$ and the center density of ${\bf E}$ is $\frac{2^r}{3^{\frac{r}{2}}}\delta({\bf L})^2$.}\\

{\bf Proof}. It follows from Theorem 4 in page 166 of \cite{CS1} directly.\\

Let $p$ be a prime number. For a lattice ${\bf L} \subset {\bf Z}^n$ we define the the linear code ${\bf C_{L,p}}$ over the finite field ${\bf F}_p$  as the image of the natural mapping ${\bf L}/p{\bf Z}^n\bigcap {\bf L} \longrightarrow {\bf F}_p^n$.\\

{\bf Theorem 2.2.} {\em Let ${\bf L} \subset {\bf Z}^n$ be a dimension $r$ lattice with minimum norm $2$. Suppose $3{\bf Z}^n \bigcap {\bf L}=3{\bf L}$. If there exists a ternary $[n,k,6]_3$ code ${\bf C}$ (over ${\bf F}_3$) which is in the ternary code ${\bf C_{L,3}}$. Then we have a lattice sphere packing of dimension $2r$ with the minimum norm $36$ and the volume $vol({\bf L})^2 \cdot 3^{\frac{5r}{2}-k}$. Its center density is $ \frac{3^k \cdot 2^r}{3^{\frac{r}{2}}}\cdot \delta({\bf L})^2$.}\\

{\bf Proof.} Let ${\bf C}'=\{({\bf c},{\bf c},{\bf c}):{\bf c} \in {\bf C}\}$ be the $[3n, k, 18]_3$ ternary code in ${\bf F}_3^n \oplus {\bf F}_3^n \oplus {\bf F}_3^n$. It is easy to check ${\bf C}'$ is in the image of the natural mapping ${\bf E}/3{\bf E}  \longrightarrow {\bf F}_3^n \oplus {\bf F}_3^n \oplus {\bf F}_3^n$. Let ${\bf T(C)}$ be the pre-image of the ternary code ${\bf C}'$.\\

We consider the lattice sphere packing of rank $2r$ in ${\bf R}^{3n}$ defined by ${\bf T(C)}$. This lattice is the union of $3^k$ translates of the lattice sphere packing $3{\bf E}$. We want to prove that the minimum distance of any two points in different translates is at least $\sqrt{36}$. Then the conclusion follows directly.\\

The residue class module $3$ in ${\bf F}_3^{3n}$  of the difference $({\bf y^1},{\bf y^2},{\bf y^3})$ (where ${\bf y^i}=(y_1^i,...,y_n^i)^{\tau}$ ) of any two vectors in different translates is of the form $({\bf c},{\bf c},{\bf c})$. Here ${\bf c}=(c_1,...c_n)^{\tau}, c_i \in {\bf F}_3=\{-1,1,0\}$, is a codeword in the ternary code ${\bf C}$. Thus the minimum norm of the difference is at least $18$. However if we look at each row component $(c_i,c_i,c_i)$ for the nonzero $c_i$. The corresponding integers $(y_i^1,y_i^2,y_i^3)$ have to satisfy $y_i^1+y_i^2+y_i^3=0$ over ${\bf Z}$. Therefore for each non-zero $c_i$ we have at least one $\pm 2$ in $(y_i^1,y_i^2,y_i^3)$. The minimum norm of $({\bf y^1},{\bf y^2},{\bf y^3})$ is at least $18+(6-3)wt({\bf c}) \geq 36$. The conclusion is proved.\\

The lattice ${\bf D}_n = \{(x_1,...,x_n): x_1+\cdots+x_n \equiv 0$ $mod$ $2 \}$ can be used in Theorem 2.2 and we get the following result.\\

{\bf Corollary 2.3.} {\em If there exists a $[n, k, 6]_3$  ternary code we have a $2n$ dimensional lattice with the volume $3^{\frac{5n-2k}{2}} \cdot 4$ and the minimum norm $36$.}\\

It is easy to check that the lattice constructed in Theorem 2.2 is always an integral even lattice.\\

Theorem 2.2 can be used to recover the lattice ${\bf E}_8$ as follows. Let ${\bf C}$ be the ternary $[8,2,6]_3$ code with the following generator matrix.\\

$$
\left(
\begin{array}{ccccccccccc}
1&-1&1&1&-1&-1&0&0\\
0&0&1&-1&1&-1&1&-1\\
\end{array}
\right)
$$

Let ${\bf L} \subset {\bf Z}^8$ be the dimension $4$ lattice defined  by ${\bf A} \cdot {\bf x}={\bf 0}$ where ${\bf A}$ is the the following integer matrix.\\

$$
\left(
\begin{array}{ccccccccccc}
1&1&0&0&0&0&0&0\\
0&0&1&0&0&1&0&0\\
0&0&0&1&1&0&0&0\\
0&0&0&0&0&0&1&1\\
\end{array}
\right)
$$

Then the volume of the lattice ${\bf L}$ is $4$ from Theorem 4 in page 166 of \cite{CS1}. The minimum norm of the lattice ${\bf L}$ is $2$. The center density of the lattice is $\delta({\bf L})=\frac{1}{2^4}$. We can check that the ternary code ${\bf C}$ is in the ternary code ${\bf C_{L,3}}$ directly. From Theorem 2.2 we get a lattice of dimension $8$ with the center density $ \frac{3^2 \cdot 2^{4}}{3^2}\cdot \frac{1}{2^{8}}=\frac{1}{16}$. Thus this is the unique densest lattice ${\bf E}_8$ of dimension 8. Then the kissing number of this lattice is $240$ (\cite{CS1}, page 123). This can be calculated from our construction directly. Actually this is another form of the tetracode construction of ${\bf E}_8$ (\cite{CS1}, page 200).\\

Similarly we have the following results.\\

{\bf Corollary 2.4.} {\em Let ${\bf L} \subset {\bf Z}^n$ be a dimension $r$ lattice with minimum norm at least $1$. Suppose $3{\bf Z}^n \bigcap {\bf L}=3{\bf L}$. If there exists a ternary $[n,k,3]_3$ code ${\bf C}$ (over ${\bf F}_3$) which is in the ternary code ${\bf C_{L,3}}$. Then we have a lattice sphere packing of dimension $2r$ with the minimum norm at least  $18$ and the volume $3^{\frac{5r}{2}-k} \cdot vol({\bf L})^2$. The center density is at least $ 3^{k-r/2} \cdot 2^{2r} \cdot \delta({\bf L})^2$.}\\

\section{The recovery of the lattice ${\bf \Lambda}_{10}$, the Kappa lattice ${\bf K}_{10}$ and the Martinet lattice ${\bf K}_{10}'$}

Let ${\bf C}$ be the ternary $[8,2,6]_3$ code with the following generator matrix.\\

$$
\left(
\begin{array}{ccccccccccc}
1&-1&1&1&-1&-1&0&0\\
0&0&1&-1&1&-1&1&-1\\
\end{array}
\right)
$$

Let ${\bf L} \subset {\bf Z}^8$ be the dimension $5$ lattice defined  by ${\bf A} \cdot {\bf x}={\bf 0}$ where ${\bf A}$ is the the following integer matrix.\\

$$
\left(
\begin{array}{ccccccccccc}
1&1&0&0&0&0&0&0\\
0&0&1&1&1&1&0&0\\
0&0&0&0&0&0&1&1\\
\end{array}
\right)
$$

Then the volume of the lattice ${\bf L}$ is $4$ from Theorem 4 in page 166 of \cite{CS1}. The minimum norm of the lattice ${\bf L}$ is $2$. The center density of the lattice is $\delta({\bf L})=\frac{1}{2^{4.5}}$. We can check that the ternary code ${\bf C}$ is in the ternary code ${\bf C_{L,3}}$ directly. Thus from Theorem 2.2 we get a lattice sphere packing ${\bf T}_{10}$ of dimension $10$ with the center density $ \frac{3^2 \cdot 2^{5}}{3^{2.5}}\cdot \frac{1}{2^9}=\frac{1}{16\sqrt{3}}$. \\

The kissing number $k_1({\bf T}_{10})=3(2+2+12)+4 \cdot 3+4 \cdot 3 \cdot 2+4 \cdot 3 \cdot 6+4 \cdot3 +4(4+4\cdot2)+4\cdot5 \cdot 6=336$ can be computed directly from the above construction. The number of second layer lattice vectors $k_2({\bf T}_{10})=48+288+432=768$ can also be computed similarly. From the recovery of the lattice ${\bf E}_8$ in the section 2 it can be proved that the lattice ${\bf T}_{10}$ is just the laminated lattice ${\bf \Lambda}_{10}$ of dimension $10$. As far as our knowledge it seems that this is the first construction of the laminated lattice of dimension $10$ from a ternary code. By equating two column coordinates we get the Kappa lattice ${\bf K}_9$ (\cite{Nebe}) with the center density $\frac{1}{16\sqrt{3}}$ and the kissing number $198$.\\

It is easy to check that the following vector ${\bf X}$ is in the real space spanned by the lattice ${\bf T}_{10}$ and the the distance of ${\bf X}$ to any lattice vector in ${\bf T}_{10}$ is at least $6$.

$$
\left(
\begin{array}{ccccccccccc}
0&0&-2&2&2&-2&0&0&\sqrt{6}\\
0&0&1&-1&-1&1&0&0&-\sqrt{6}\\
0&0&1&-1&-1&1&0&0&0\\
\end{array}
\right)
$$
The lattice generated by  ${\bf X}$ and the lattice ${\bf T}_{10}$ is a lattice with the center density $\frac{1}{32}$ and the kissing number $336+2 \cdot (4 \cdot 3 \cdot 3+15)=438$. This construction is actually based on a ternary $[8, 3, 4]_3$ code. This is the lattice ${\bf \Lambda}_{11}^{max}$ (\cite{Nebe}).\\

The Kappa lattice ${\bf K}_{10}$ can be recovered easily as follows.  Let ${\bf C}$  be the trivial $[6, 1, 6]_3$ ternary code generated by $(1,1,1,-1,-1,-1)$ and ${\bf L} \subset {\bf Z}^6$ defined by $x_1+\cdots+x_6=0$. Then from Theorem 2.2 we get a $10$ dimensional lattice with the center density $\frac{1}{18\sqrt{3}}$. Its kissing number is $90+2(3+ 9 \cdot 4 + 9\cdot 2 \cdot3)=276$. This $10$ dimensional lattice is the Kappa lattice ${\bf K}_{10}$. If ${\bf C}$ is the trivial $[6, 1, 6]_3$ ternary code generated by $(1,1,1,1,1,1)$ and ${\bf L} \subset {\bf Z}^6$ is the $5$ dimensional lattice defined by $x_1+\cdots+x_6=0$. From Theorem 2.2 we get a $10$ dimensional lattice with the same center density $\frac{1}{18\sqrt{3}}$. Its kissing number is $90+2 \cdot 15 \cdot 6=270$. This is the Martinet lattice ${\bf K}_{10}'$ (\cite{Mar}, page 290).

\section{The recovery of the Coxeter-Todd  lattice}

Let ${\bf C}$ be the ternary $[9,3,6]_3$ code with the following generator matrix.\\

$$
\left(
\begin{array}{ccccccccccc}
1&0&-1&0&1&-1&-1&0&1\\
0&0&0&1&1&1&-1&-1&-1\\
0&1&-1&0&-1&1&-1&1&0\\
\end{array}
\right)
$$

Let ${\bf L} \subset {\bf Z}^9$ be the dimension $6$ lattice defined  by ${\bf A} \cdot {\bf x}={\bf 0}$ where ${\bf A}$ is the the following integer matrix.\\

$$
\left(
\begin{array}{ccccccccccc}
1&1&1&0&0&0&0&0&0\\
0&0&0&1&1&1&0&0&0\\
0&0&0&0&0&0&1&1&1\\
\end{array}
\right)
$$

Then the volume of the lattice ${\bf L}$ is $\sqrt{27}$ from Theorem 4 in page 166 of \cite{CS1}. The minimum norm of the lattice ${\bf L}$ is $2$. The center density of the lattice is $\delta({\bf L})=\frac{1}{2^3\sqrt{27}}$. We can check that the ternary code ${\bf C}$ is in the ternary code ${\bf C_{L,3}}$ directly. Thus from Theorem 2.3 we get a lattice  ${\bf K}_{12}$ of dimension $12$ with the center density $\frac{3^3 \cdot 2^6}{3^3}\cdot \frac{1}{2^6\cdot 27}=\frac{1}{27}$. This is just the Coxeter-Todd lattice of dimension $12$.\\

One base of the above lattice ${\bf K}_{12}$ is as follows. Three rows are listed in a same row.\\
$[-2,0,2,0,-2,2,2,0,-2,1,0,-1, 0,1,-1, -1,0,1, 1,0,-1,0,1,-1,-1,0,1]$;\\
$[0,-2,2,0,2-2,2,-2,0,  0,1,-1, 0,-1,1,-1,1,0,  0,1,-1, 0,-1,1,-1,1,0]$;\\
$[0,0,0,-2,1,1,2,-1,-1, 0,0,0,1,-2,1,-1,2,-1, 0,0,0, 1,1,-2,-1,-1,2]$;\\
$[3,-3,0,0,0,0,0,0,0,-3,3,0,0,0,0,0,0, 0, 0,0,0,0,0,0,0,0,0]$;\\
$[3,-3,0,0,0,0,0,0,0, 0,0,0,0,0,0,0,0,0,-3,3,0,0,0,0,0,0,0]$;\\
$[3,0,-3,0,0,0,0,0,0,-3,0,3,0,0,0,0,0, 0,0,0,0,0,0,0,0,0,0]$;\\
$[3,0, -3,0,0,0,0,0,0, 0,0,0,0,0,0,0,0,0,-3,0,3,0,0,0,0,0,0]$;\\
$[0,0,0,0,3,-3, 0,0,0, 0, 0,0, 0,-3,3,0,0,0,0, 0,0, 0,0,0,0,0,0]$;\\
$[0,0,0,0,3,-3, 0,0,0, 0, 0,0, 0,0,,0, 0,0,0, 0,0,0, 0,-3,3, 0,0,0]$;\\
$[0,0,0, 3,0,-3,0,0,0, 0,0,0,-3,0,3,0,0,0 ,0,0, 0,0,0,0,0,0,0]$;\\
$[0,0,0,0,0,0, 3,-3,0, 0,0,0, 0,0,0, -3,3,0, 0,0,0,0,0,0,0,0,0]$;\\
$[0,0,0, 0,0,0, 3,0,-3, 0,0,0, 0,0,0 -3,0,3,  0,0,0, 0,0,0, 0,0,0]$.\\
With the help of the Magma the Gram matrix is as follows.\\
$$
\left(
\begin{array}{cccccccccccccccccccccc}
36&0&   9&  -9&  -9 &-18& -18 &-18 &-18 & -9 &  9  &18\\
0 & 36&   9&   9 &  9  &-9  &-9  &18 & 18&   9  &18  & 9\\
 9 &  9  &36  & 0  & 0  & 0&   0  & 9  &-9  &-9 & 18 &  9\\
-9 &  9 &  0&  36 & 18  &18  & 9  & 0&   0 &  0  & 0 &  0\\
-9&   9&   0&  18&  36&   9&  18&   0&   0&   0&   0&   0\\
-18&  -9&   0&  18&   9&  36&  18&   0&   0&   0&   0&   0\\
-18 & -9 &  0 &  9&  18&  18&  36&   0&   0&   0&   0&   0\\
-18 & 18 &  9&   0&   0&   0&   0&  36&  18&  18&   0&   0\\
-18  &18 & -9   &0 &  0 &  0 &  0 & 18&  36 &  9&   0&   0\\
 -9 &  9  &-9&   0   &0 &  0 &  0 & 18  & 9  &36  & 0  & 0\\
 9 & 18  &18 &  0   &0  & 0 &  0  & 0  & 0  & 0 & 36 & 18\\
18 &  9 &  9 &  0 &  0 &  0  & 0 &  0  & 0  & 0  &18 & 36\\
\end{array}
\right)
$$

 We can get an isometry to the Coxeter-Todd lattice ${\bf K}_{12}$ by Magma.\\

If we use the $5$ dimensional lattice ${\bf L} \subset {\bf Z}^8$ defined by ${\bf A} \cdot {\bf x}=0$ where ${\bf A}$ is the following matrix

$$
\left(
\begin{array}{ccccccccccc}
1&1&0&0&0&0&0&0\\
0&0&1&1&1&0&0&0\\
0&0&0&0&0&1&1&1\\
\end{array}
\right)
$$
and the $[8, 2, 6]_3$ ternary code ${\bf C}$ with the following generator matrix
$$
\left(
\begin{array}{ccccccccccc}
1&-1&1&-1&0&1&-1&0\\
1&-1&0&1&-1&0&1&-1\\
\end{array}
\right)
$$
we get the lattice ${\bf K}_{10}$ with the center density $\frac{1}{18\sqrt{3}}$ and the kissing number $276$. The distance of the following vector ${\bf X}$ in the space spanned by the lattice ${\bf K}_{10}$ to any lattice  vector in ${\bf K}_{10}$ is at least $\sqrt{27+9}=6$.

$$
\left(
\begin{array}{ccccccccccccccccc}
0&0&-1&-1&2&1&1&-2&0\\
\frac{3}{2}&-\frac{3}{2}&-1&2&-1&1&-2&1&\sqrt{\frac{9}{2}}\\
-\frac{3}{2}&\frac{3}{2}&2&-1&-1&-2&1&1&-\sqrt{\frac{9}{2}}\\
\end{array}
\right)
$$

The lattice generated by ${\bf K}_{10}$ and this ${\bf X}$ ($m{\bf X}+{\bf K}_{10}$) is  a lattice ${\bf K}_{11}$ with the center density $\frac{1}{18\sqrt{3}}$ and the kissing number $432$.\\

Let ${\bf Y}$ be the following vector.\\

$$
\left(
\begin{array}{ccccccccccccccccccc}
\frac{3}{2}&-\frac{3}{2}&-1&-1&2&1&1&-2&\sqrt{\frac{9}{2}}\\
-\frac{3}{2}&\frac{3}{2}&-1&2&-1&1&-2&1&-\sqrt{\frac{9}{2}}\\
0&0&2&-1&-1&-2&1&1&0\\
\end{array}
\right)
$$

Here we should note that the following generator matrix corresponds to a $[8, 3, 3]_3$ ternary code ${\bf C}'$.\\

$$
\left(
\begin{array}{ccccccccccc}
1&-1&1&-1&0&1&-1&0\\
1&-1&0&1&-1&0&1&-1\\
0&0&-1&-1&-1&1&1&1\\
\end{array}
\right)
$$

In the teranry code ${\bf C}'$, there are four Hamming weight $3$ codewords, $8(in {\bf C})+12(outside {\bf C})+2(3rd-row)$ Hamming weight $6$ codewords.\\

The lattice generated by ${\bf X}, {\bf Y}$ and ${\bf K}_{10}$ is the $12$ dimensional Coxeter-Todd lattice ${\bf K}_{12}$ again. Its volume is $18 \sqrt{3} \cdot \sqrt{9} \cdot \sqrt{\frac{3}{2} \cdot \frac{9}{2}}=3^{15}$ and the minimum norm is $36$. Its center density is $\frac{1}{27}$ and its kissing number is $756=276+6(3 \cdot3 \cdot3 +3\cdot 3 \cdot3 +2 \cdot 2 \cdot 6)+6 \cdot 2$. Here $276$ is the kissing number of the lattice ${\bf K}_{10}$, the two $3 \cdot 3 \cdot 3$'s correspond to the half of the $12$ Hamming weight $6$ outside codewords. Here we should note that the 3rd row lead to no norm $36$ lattice vectors. The last $2 \cdot 6 \cdot 2$ corresponds to the half of the four Hamming weight $3$ codewords in ${\bf C}'$. The construction of ${\bf K}_{12}$ is similar to the construction of ${\bf \Lambda}_{14}$ in section 8.\\

The gram matrix of the lattice ${\bf K}_{12}$ is as follows.\\
$$
\left(
\begin{array}{ccccccccccccccccccccc}
36&   0& -18& -18& -18&  9& -18&  9& -18& -9&   0& -9\\
 0 & 36 &-18& -18&  9& -18&  9& -18&  9& -9&   0& -9\\
-18 &-18&  36&  18&   0&   0&   0&   0&   0&   0& -9&  18\\
-18& -18&  18& 36&   0&   0&   0&   0&   0&   0&  9&  9\\
-18&  9 &  0&   0&  36& -18&  18& -9&   0&   0&  9&  9\\
9& -18&   0&   0& -18&  36& -9 & 18&   0&   0& -18& -18\\
-18&  9&   0&   0&  18& -9&  36& -18&   0&   0& -9& -9\\
 9& -18&   0&   0& -9&  18& -18&  36&   0&   0& -9& -9\\
-18 & 9 &  0&   0&   0&   0&   0&   0&  36&  18& -9& -9\\
-9& -9 &  0 &  0&   0&   0&   0&   0&  18&  36&  9&  9\\
0   &0 &-9 & 9&  9& -18& -9& -9& -9&  9&  54&  27\\
-9 &-9&  18&  9&  9& -18& -9& -9& -9&  9&  27& 54\\
\end{array}
\right)
$$

By Magma we can get an isometry to the Coxeter-Todd lattice ${\bf K}_{12}$. \\

\section{The recovery of the Barnes-Wall lattice}

Let ${\bf C}$ be the following $[8, 4, 3]_3$ ternary code with the following generator matrix.\\

$$
\left(
\begin{array}{cccccccccccccccc}
1&1&1&0&0&0&0&0\\
0&1&-1&1&0&0&0&0\\
0&0&0&0&1&1&1&0\\
0&0&0&0&0&1&-1&1\\
\end{array}
\right)
$$

There are $16$ Hamming weight $3$ codewords and $64$ Hamming weight $6$ codewords in ${\bf C}$.\\

Let $L \subset ({\bf Z}^8 \oplus {\bf Z}^8 \oplus {\bf Z}^8)$ be a $16$ dimensional lattice defined by the column conditions $x_{1j}+x_{2j}+x_{3j}=0$, where $j=1,...,8$ and $x_{ij}, i=1,2,3,j=1,...,8$'s are the coordinates of the lattice ${\bf Z}^8 \oplus {\bf Z}^8 \oplus {\bf Z}^8$. From Theorem 2.3 we get a $16$ dimensional lattice ${\bf T'}_{16}$ with the minimum norm $18$ and the volume $\frac{3^{4} \cdot 1^2 \cdot 3^{16}}{3^4}=3^{16}$. Its center density is $\frac{1}{256}$. The following vectors ${\bf 0}$, ${\bf X}_1,...,{\bf X}_7$ \\

$$
\left(
\begin{array}{cccccccccccccccc}
\frac{3}{2}&\frac{3}{2}&0&0&\frac{3}{2}&\frac{3}{2}&0&0\\
-\frac{3}{2}&-\frac{3}{2}&0&0&-\frac{3}{2}&-\frac{3}{2}&0&0\\
0&0&0&0&0&0&0&0\\
\end{array}
\right)
$$

$$
\left(
\begin{array}{cccccccccccccccc}
0&\frac{3}{2}&\frac{3}{2}&0&0&\frac{3}{2}&\frac{3}{2}&0\\
0&-\frac{3}{2}&-\frac{3}{2}&0&0&-\frac{3}{2}&-\frac{3}{2}&0\\
0&0&0&0&0&0&0&0\\
\end{array}
\right)
$$

$$
\left(
\begin{array}{cccccccccccccccc}
0&0&\frac{3}{2}&\frac{3}{2}&0&\frac{3}{2}&0&\frac{3}{2}\\
0&0&-\frac{3}{2}&-\frac{3}{2}&0&-\frac{3}{2}&0&-\frac{3}{2}\\
0&0&0&0&0&0&0&0\\
\end{array}
\right)
$$

$$
\left(
\begin{array}{cccccccccccccccc}
\frac{3}{2}&0&\frac{3}{2}&0&\frac{3}{2}&\frac{3}{2}&\frac{3}{2}&\frac{3}{2}\\
-\frac{3}{2}&0&-\frac{3}{2}&0&-\frac{3}{2}&-\frac{3}{2}&-\frac{3}{2}&-\frac{3}{2}\\
0&0&0&0&0&0&0&0\\
\end{array}
\right)
$$

$$
\left(
\begin{array}{cccccccccccccccc}
\frac{3}{2}&\frac{3}{2}&\frac{3}{2}&\frac{3}{2}&\frac{3}{2}&0&0&\frac{3}{2}\\
-\frac{3}{2}&-\frac{3}{2}&-\frac{3}{2}&-\frac{3}{2}&-\frac{3}{2}&0&0&-\frac{3}{2}\\
0&0&0&0&0&0&0&0\\
\end{array}
\right)
$$

$$
\left(
\begin{array}{cccccccccccccccc}
0&\frac{3}{2}&0&\frac{3}{2}&0&\frac{3}{2}&\frac{3}{2}&0\\
0&-\frac{3}{2}&0&-\frac{3}{2}&0&-\frac{3}{2}&-\frac{3}{2}&0\\
0&0&0&0&0&0&0&0\\
\end{array}
\right)
$$

$$
\left(
\begin{array}{cccccccccccccccc}
\frac{3}{2}&0&0&\frac{3}{2}&\frac{3}{2}&0&\frac{3}{2}&0\\
-\frac{3}{2}&0&0&-\frac{3}{2}&-\frac{3}{2}&0&-\frac{3}{2}&0\\
0&0&0&0&0&0&0&0\\
\end{array}
\right)
$$
and the following vector ${\bf Y}$ \\

$$
\left(
\begin{array}{cccccccccccccccc}
0&0&0&0&0&0&0&0\\
-\frac{3}{2}&-\frac{3}{2}&-\frac{3}{2}&-\frac{3}{2}&-\frac{3}{2}&-\frac{3}{2}&-\frac{3}{2}&-\frac{3}{2}\\
\frac{3}{2}&\frac{3}{2}&\frac{3}{2}&\frac{3}{2}&\frac{3}{2}&\frac{3}{2}&\frac{3}{2}&\frac{3}{2}\\
\end{array}
\right)
$$
are needed. It is easy to check that the above $8$ vectors are in the real space ${\bf R}^{16}$ spanned by the lattice ${\bf T'}_{16}$.\\

We consider the $16$ vectors of the form ${\bf 0},{\bf X}_1,{\bf X}_2,{\bf X}_3$,${\bf X}_4,{\bf X}_5$, ${\bf X}_6,{\bf X}_7$,${\bf Y}, {\bf Y}-{\bf X}_1,{\bf Y}-{\bf X}_2,{\bf Y}-{\bf X}_3$,${\bf Y}-{\bf X}_4,{\bf Y}-{\bf X}_5,{\bf Y}-{\bf X}_6,{\bf Y}-{\bf X}_7$. The union of the sixteen translates of the above lattice ${\bf T'}_{16}$ is a lattice ${\bf T}_{16}$ with the volume $\frac{3^{16}}{16}$. Here it should be noted that $2{\bf X}_i, 2{\bf Y} \in {\bf L} \subset {\bf T'}_{16}$.\\

{\bf Theorem 5.1.} {\em The minimum norm of the lattice ${\bf T}_{16}$ is $18$.}\\

{\bf Proof.} It is clear that the Euclid norms of the differences of the vectors of the form $3 {\bf v}$ to any above sixteen vectors are at least $\frac{9}{2} \cdot 4=18$. For any vector of the form ${\bf V}=({\bf v}_1,{\bf v}_2,{\bf v}_3)$ where ${\bf v}_1 \equiv {\bf v}_2 \equiv {\bf v}_3 \equiv {\bf c} \in {\bf C}$ a non-zero codeword with Hamming weight $3$ or $6$ in the above $[8, 4, 3]_3$ ternary code, the Euclid norm of the differences of this ${\bf V}$ to the seven vectors ${\bf 0},{\bf X}_1,...,{\bf X}_7$ are at least $\frac{3}{2} \cdot 2 +\frac{9}{2} \cdot 2=\frac{3}{2} \cdot 2 +6 +\frac{3}{2} \cdot2 +6=\frac{9}{2}+3 \cdot \frac{3}{2}+\frac{9}{2} \cdot 2=18$.  The Euclid norm of the differences of this ${\bf V}$ to the nine vectors ${\bf T}_7, {\bf Y},{\bf Y}-{\bf X}_1,...,{\bf Y}-{\bf X}_7$ are at least $\frac{3}{2} \cdot 3 + \frac{9}{2} \cdot 5=27$(for Hamming weight $3$ codewords)  or at least $\frac{3}{2} \cdot 6 +\frac{9}{2} \cdot 2==18$(for Hamming weight $6$ codewords). This can be easily checked from the analysis in the section 2.\\

The following is a base of the lattice ${\bf T}_{16}$. We can get an isometry to the Barnes-Wall lattice ${\bf \Lambda}_{16}$ by Magma.\\
$[3, 0, 0, 0, 0, 0, 0, 0, 0, 0, 0, 0, 0, 0, 0, 0, -3, 0, 0, 0, 0, 0, 0, 0]$\\
$[0, 0, 0, 0, 0, 0, 0, 0, 3, 0, 0, 0, 0, 0, 0, 0, -3, 0, 0, 0, 0, 0, 0, 0]$\\
$[1, 2, 0, -2, 0, 0, 0, 0, 1, -1, 0, 1, 0, 0, 0, 0, -2, -1, 0, 1, 0, 0, 0, 0]$\\
$[1, -1, 0, 1, 0, 0, 0, 0, 1, -1, 0, 1, 0, 0, 0, 0, -2, 2, 0, -2, 0, 0, 0, 0]$\\
$[-1, 1, 0, 2, 0, 0, 0, 0, -1, 1, 0, -1, 0, 0, 0, 0, 2, -2, 0, -1, 0, 0, 0, 0]$\\
$[1, -1, 0, -2, 0, 0, 0, 0, 1, 2, 0, 1, 0, 0, 0, 0, -2, -1, 0, 1, 0, 0, 0, 0]$\\
$[1, 1, -2, 0, 0, 0, 0, 0, 1, 1, 1, 0, 0, 0, 0, 0, -2, -2, 1, 0, 0, 0, 0, 0]$\\
$[1, 1, 1, 0, 0, 0, 0, 0, 1, 1, 1, 0, 0, 0, 0, 0, -2, -2, -2, 0, 0, 0, 0, 0]$\\
$[\frac{3}{2}, \frac{3}{2}, 0, 0, \frac{3}{2}, \frac{3}{2}, 0, 0, -\frac{3}{2}, -\frac{3}{2}, 0, 0, -\frac{3}{2}, -\frac{3}{2}, 0, 0, 0, 0, 0, 0, 0, 0, 0, 0]$\\
$[-\frac{1}{2}, 1, -\frac{1}{2}, 0, \frac{1}{2}, \frac{3}{2}, -\frac{1}{2}, -\frac{1}{2}, -\frac{1}{2}, 1, -\frac{1}{2}, 0, \frac{1}{2}, -\frac{3}{2}, -\frac{1}{2}, -\frac{1}{2}, 1, -2, 1, 0, -1, 0, 1, 1]$\\
$[2, \frac{1}{2}, \frac{1}{2}, 0, -2, 0, \frac{1}{2}, \frac{1}{2}, -1, \frac{1}{2}, \frac{1}{2}, 0, 1, 0, \frac{1}{2}, \frac{1}{2}, -1, -1, -1, 0, 1, 0, -1, -1]$\\
$[2, \frac{1}{2}, \frac{1}{2}, 0, 0, 2, -\frac{1}{2}, \frac{1}{2}, -1, \frac{1}{2}, \frac{1}{2}, 0, 0, -1, -\frac{1}{2}, \frac{1}{2}, -1, -1, -1, 0, 0, -1, 1, -1]$\\
$[-\frac{3}{2}, 0, 0,  \frac{3}{2}, -\frac{1}{2}, 1, -\frac{1}{2}, 0,  \frac{3}{2}, 0, 0, -\frac{3}{2}, -\frac{1}{2}, -2, -\frac{1}{2}, 0, 0, 0, 0, 0, 1, 1, 1, 0]$\\
$[\frac{1}{2}, \frac{1}{2}, \frac{1}{2}, -\frac{3}{2}, -\frac{3}{2}, -1, 1, \frac{1}{2}, -1, -1, -1, 0, 0, \frac{1}{2}, -\frac{1}{2}, -1, \frac{1}{2}, \frac{1}{2}, \frac{1}{2}, \frac{3}{2}, \frac{3}{2},  \frac{1}{2}, -\frac{1}{2}, \frac{1}{2}]$\\
$[1, 0, -1, -1, 1, 1, 1, 0, -\frac{1}{2}, -\frac{3}{2}, \frac{1}{2}, \frac{1}{2}, -\frac{1}{2}, -\frac{1}{2}, -\frac{1}{2}, \frac{3}{2}, -\frac{1}{2}, \frac{3}{2}, \frac{1}{2}, \frac{1}{2},  -\frac{1}{2}, -\frac{1}{2}, -\frac{1}{2}, -\frac{3}{2}]$\\
$[0, -\frac{1}{2}, \frac{1}{2}, 1, 1, 0, \frac{1}{2}, \frac{1}{2}, \frac{3}{2}, 1, -1, -\frac{1}{2}, -\frac{1}{2}, -\frac{3}{2}, -1, -1, -\frac{3}{2}, -\frac{1}{2}, \frac{1}{2}, -\frac{1}{2}, - \frac{1}{2}, \frac{3}{2},  \frac{1}{2}, \frac{1}{2}]$\\

From some normal vectors the $15$ dimensional laminated lattice ${\bf \Lambda}_{15}$ can be constructed similarly.\\

\section{The recovery of the Leech lattice}

The $[12, 6, 6]_3$ ternary Golay code is defined by the following generator matrix (\cite{CS1}, page 85).\\

$$
\left(
\begin{array}{cccccccccccccc}
1&0&0&0&0&0&0&1&1&1&1&1\\
0&1&0&0&0&0&-1&0&1&-1&-1&1\\
0&0&1&0&0&0&-1&1&0&1&-1&-1\\
0&0&0&1&0&0&-1&-1&1&0&1&-1\\
0&0&0&0&1&0&-1&-1&-1&1&0&1\\
0&0&0&0&0&1&-1&1&-1&-1&1&0\\
\end{array}
\right)
$$
It is a self-dual code with weight distribution $A_0=1, A_6=264, A_9=440, A_{12}=24$. Here $A_i$ is the number of codewords of Hamming weight $i$. The $24$ weight $12$ codewords are of the form $\pm(1,...,1)$ (two codewords) or $\pm(1,...,1,-1,...,-1)$ ($22$ codewords have six $1$'s and six $-1$'s). The $440$ weight $9$ codewords are of the form $\pm (1^6,(-1)^3,0^3)$ ($220$ codewords have six $1$'s and three $-1$'s and $220$ codewords have six $-1$'s and three $1$'s). The $264$ weight $6$ codewords are of the form $\pm(1^6)$ ($22$ codewords have six $1$'s and $22$ codewords have six $-1$'s) or $(1^3,(-1)^3)$ ($220$ codewords have three $1$'s and three $-1$'s). From Corollary 2.3 we have a $24$ dimensional lattice ${\bf T}_{24}'$ with the volume $3^{24} \cdot 4$ and the minimum norm $36$. \\

Let ${\bf x}=(\frac{3}{2},\frac{3}{2},\frac{3}{2},\frac{3}{2},\frac{3}{2},\frac{3}{2},\frac{3}{2},-\frac{3}{2},-\frac{3}{2},-\frac{3}{2},-\frac{3}{2},-\frac{3}{2})$, and ${\bf X}_1=({\bf x}, -{\bf x},{\bf 0}), {\bf X}_2\\=(-{\bf x}, {\bf 0},{\bf x}), {\bf X}_3=({\bf 0}, -{\bf x},{\bf x})$. We want to prove that the union of the four tanslates ${\bf T}_{24}', {\bf X}_1+{\bf T}_{24}', {\bf X}_2+{\bf T}_{24}', {\bf X}_3+{\bf T}_{24}'$ of ${\bf T}_{24}'$ is a lattice with the volume $3^{24}$ and the minimum norm $36$. It is easy to check this is a lattice ${\bf T}_{24}$ with volume $3^{24}$ since $2{\bf X}_i \in {\bf T}_{24}'$. Here we can verify that the vector ${\bf x}$ can be replaced by any vector with odd numbers of $\frac{3}{2}$'s and $-\frac{3}{2}$'s and the resulted lattice is the same. \\

The only remaining point is to determine the minimum norm of the lattice ${\bf T}_{24}$. Here we check that the Euclid norm of each vector in the translate ${\bf X}_1+{\bf T}_{24}'$ is at least $36$. It is clear that the Euclid norm of vectors in ${\bf X}_1+3 {\bf E}$ is at least $54$. For any vector ${\bf V} \in {\bf T}_{24}'$ with its residue class equal to a Hamming weight $6$ codeword in the $[12, 6, 6]_3$ extended ternary Golay code, the Euclid norm of ${\bf X}_1-{\bf V}$ is at least $6 \cdot (\frac{9}{4}+\frac{9}{4})+ 6 \cdot (\frac{1}{4}+\frac{1}{4}+1)=36$. When ${\bf V}=({\bf v}_1,{\bf v}_2,{\bf v}_3)$, where the Hamming weight of the residue classes of all ${\bf v}_i$'s is $9$, there have to be at least one $\pm2$ in each ${\bf v}_i$ since ${\bf v}_i$ is in ${\bf D}_{12}$ (see Corollary 2.3 and the construction in Theorem 2.2) or at least one $\pm3$ at the last column, the Euclid norm of ${\bf X}_1-{\bf V}$ is at least $3 \cdot (\frac{9}{4}+\frac{9}{4})+ 9 \cdot (\frac{1}{4}+\frac{1}{4}+1)+9=36$. When ${\bf V}=({\bf v}_1,{\bf v}_2,{\bf v}_3)$, where the Hamming weight of the residue classes of all ${\bf v}_i$'s is $12$, there have to be at least one $\pm\frac{5}{2}$ and $\pm \frac{7}{2}$ in each ${\bf X}_1-{\bf V}$ since there are odd numbers of $\pm \frac{3}{2}$'s  in ${\bf x}$, or at least one pair of $\pm(\frac{1}{2},-\frac{5}{2},2)$ or $\pm(-\frac{5}{2},\frac{1}{2},2)$, the Euclid norm of ${\bf X}_1-{\bf V}$ is at least $12 \cdot (\frac{1}{4}+\frac{1}{4}+1)+9+9=12 \cdot (\frac{1}{4}+\frac{1}{4}+1)+18=36$. Similar argument is also valid for other translates of ${\bf T}_{24}'$. Thus the lattice ${\bf T}_{24}$ is of the  volume $3^{24}$ and the minimum norm $36$. The center density of the lattice ${\bf T}_{24}$ is $1$.\\

It is clear that $\frac{1}{3}{\bf T}_{24}$ is an unimodualr even lattice with the minimum norm $4$. From the characterization of the Leech lattice (Chapter 12 of \cite{CS1}) this is the Leech lattice. This can also be proved simply from the celebrated theorem of \cite{CoKu}.\\

It is well-known that the Leech lattice can be constructed from the extended $[12, 6, 6]_3$ ternary Golay code and a lattice over the ring of the Eisenstein integers (\cite{CS1}, page 200). We do not know the relation between our this construction and that number-theoretic construction. In both constructions the ternary Golay code is a main ingredient.\\

 From some normal vectors, the $23$ dimensional laminated lattice ${\bf \Lambda}_{23}$ can be constructed similarly.\\

\section{The recovery of the laminated lattices ${\bf \Lambda}_{22}$ and ${\bf \Lambda}_{26}$}

Let ${\bf L} \subset {\bf Z}^{12}$ be a $11$ dimensional lattice of the form ${\bf D}_{10} \oplus {\bf A}_1$, where ${\bf A}_1$ is defined by $x_1-x_2=0$. It is easy to verify that each row of the following generator matrix of the ternary $[12, 5, 6]_3$ code  is in ${\bf C_{L,3}}$. The weight distribution of this $[12, 5, 6]_3$ ternary code is $A_0=1,A_6=90,A_{9}=140, A_{12}=12$.\\
$$
\left(
\begin{array}{cccccccccccccc}
1&1&0&0&0&0&-1&1&-1&0&0&-1\\
0&0&1&0&0&0&-1&1&0&1&-1&-1\\
0&0&0&1&0&0&-1&-1&1&0&1&-1\\
0&0&0&0&1&0&-1&-1&-1&1&0&1\\
0&0&0&0&0&1&-1&1&-1&-1&1&0\\
\end{array}
\right)
$$
From Theorem 2.2 we get a lattice ${\bf T}_{22}'$ with the center density $\frac{1}{8\sqrt{3}}$. Let ${\bf x}=(\frac{3}{2},\frac{3}{2},\frac{3}{2},\frac{3}{2},\frac{3}{2},\frac{3}{2},\frac{3}{2},\frac{3}{2},\frac{3}{2},
-\frac{3}{2},-\frac{3}{2},-\frac{3}{2})$, and ${\bf X}_1=({\bf x}, -{\bf x},{\bf 0}), {\bf X}_2=(-{\bf x}, {\bf 0},{\bf x}), {\bf X}_3\\=({\bf 0}, -{\bf x},{\bf x})$. As argued in the previous section the union of the four translates ${\bf T}_{22}', {\bf X}_1+{\bf T}_{22}', {\bf X}_2+{\bf T}_{22}', {\bf X}_3+{\bf T}_{22}'$ of ${\bf T}_{22}'$ is a lattice ${\bf T}_{22}$ with the volume $2\cdot3^{22.5}$ and the minimum norm is $36$. The lattice $\frac{1}{3}{\bf T}_{22}$ is an integral even lattice with the determinant $12$ and the minimum norm $4$. Since ${\bf T}_{22}$ is a cross section defined by two minimum norm vectors in ${\bf T}_{24}$ with the inner product $2$ it is just the laminated lattice ${\bf \Lambda}_{22}$, since the automorphism group of the Leech lattice is transitive on the set of pairs of minimum norm lattice vectors with the fixed inner product (\cite{CS1}, Chapter 10).\\

It can be verified similarly as the previous section that the kissing number of the lattice ${\bf T}_{22}$ is $k_1({\bf T}_{22})=k_1({\bf T}_{22}')+3(20 \cdot 2+ 9 \cdot 10 \cdot 12+ 120 \cdot 12$$+ 2 \cdot 10 \cdot 12+20 \cdot 4 \cdot 9+ 120 \cdot 7 \cdot 8+ 30 \cdot 2^5 +60 \cdot 2^4)=k_1({\bf T}_{22}')+ 3\cdot 12160$, where $k_1({\bf T}_{22}')= 182 \cdot3 + 60 \cdot 183 +30 \cdot 63=13416$. Thus $k_1({\bf T}_{22})=13416+3 \cdot 12160=49896$.\\

It is well-known that there is an unique integral even lattice of dimension $26$ with the determinant $3$ and the minimum norm $4$ (\cite{CS1,Ba}). This lattice is also one of the known densest (laminated) lattice of dimension $26$ (\cite{CS1}, page xix). This lattice can be constructed naturally from our main result Theorem 2.2 and the ternary $[13, 6, 6]_3$ code.\\

We consider the $[13, 6, 6]_3$ ternary code with the following generator matrix. It contains the $[12, 6, 6]_3$ ternary Golay code as a sub-code.\\

$$
\left(
\begin{array}{cccccccccccccc}
1&0&0&0&0&0&0&1&1&1&1&1&0\\
0&1&0&0&0&0&-1&0&1&-1&-1&1&0\\
0&0&1&0&0&0&-1&1&0&1&-1&-1&0\\
0&0&0&1&0&0&-1&-1&1&0&1&-1&0\\
0&0&0&0&1&0&-1&-1&-1&1&0&1&0\\
0&0&0&0&0&1&-1&1&-1&-1&1&0&0\\
\end{array}
\right)
$$
From Corollary 2.3 we have a $26$ dimensional lattice ${\bf T}_{26}'$ with the volume $3^{26.5} \cdot 4$ and the minimum norm $36$. Let ${\bf x}=(\frac{3}{2},\frac{3}{2},\frac{3}{2},\frac{3}{2},\frac{3}{2},-\frac{3}{2},-\frac{3}{2},-\frac{3}{2},-\frac{3}{2},-\frac{3}{2},\\-\frac{3}{2},-\frac{3}{2},0)$, and ${\bf X}_1=({\bf x}, -{\bf x},{\bf 0}), {\bf X}_2=(-{\bf x}, {\bf 0},{\bf x}), {\bf X}_3=({\bf 0}, -{\bf x},{\bf x})$. The union of the four translates ${\bf T}_{26}', {\bf X}_1+{\bf T}_{26}', {\bf X}_2+{\bf T}_{26}', {\bf X}_3+{\bf T}_{26}'$ of ${\bf T}_{26}'$ is a lattice ${\bf T}_{26}$ with the volume $3^{26.5}$, since $2{\bf X}_i \in {\bf T}_{26}'$. From a similar argument as in the section 6 the minimum Euclid norm of the lattice ${\bf T}_{26}$ is $36$.\\

Since there are $\frac{13 \cdot 12}{2}\cdot 4=312$ Euclid norm $2$ lattice vectors in the lattice ${\bf D}_{13}$, argued as the case of lattice ${\bf T}_{24}$, $k_1({\bf T}_{26}')=3 \cdot 312 + 264 \cdot 183=49248$. Similarly as the case of the lattice ${\bf T}_{24}$ and from the above argument the kissing number $k_1({\bf T}_{26})=49248 + 3 \cdot (264 \cdot 32 +440\cdot 72+440 \cdot 12 +24 \cdot 12 \cdot 11+ 2 \cdot 12 \cdot 24 +24\cdot 2)= =49248+3 \cdot 49200=196848$. This is the same as the highest known kissing number in ${\bf R}^{26}$ which is attained by one laminated lattice of dimension $26$. It is easy to check that the lattice $\frac{1}{3}{\bf T}_{26}$ is an integral even lattice with the determinant $3$ and the minimum norm $4$. It is not hard to verify that the lattice ${\bf T}_{26}$ is the same as the laminated lattice ${\bf \Lambda}_{26}$.\\

\section{The recovery of the laminated lattice ${\bf \Lambda}_{14}$}

Let ${\bf L} \subset {\bf Z}^8$ be the dimension $6$ lattice defined  by ${\bf A} \cdot {\bf x}={\bf 0}$ where ${\bf A}$ is the the following integer matrix.\\

$$
\left(
\begin{array}{ccccccccccc}
1&1&0&0&0&0&-1&-1\\
0&0&1&-1&-1&1&0&0\\
\end{array}
\right)
$$

Let ${\bf C} \subset {\bf C}'$ be the ternary $[8,2,6]_3$ code and $[8, 3, 5]_3$ code with the following generator matrices.
There are $2$ Hamming weight $8$ codewords and $16$ Hamming weight $5$ codewords in ${\bf C'}$.\\

$$
\left(
\begin{array}{ccccccccccc}
1&-1&1&1&-1&-1&0&0\\
0&0&1&-1&1&-1&1&-1\\
\end{array}
\right)
$$

$$
\left(
\begin{array}{ccccccccccc}
1&-1&1&1&-1&-1&0&0\\
0&0&1&-1&1&-1&1&-1\\
1&1&1&1&1&1&1&1\\
\end{array}
\right)
$$

It is clear that ${\bf C} \subset {\bf C}'$ are two ternary codes in the ternary code ${\bf C_{L,3}}$. From Theorem 2.2 we have a $12$ dimensional lattice ${\bf W}_{12}$ with the volume $3^{13} \cdot 16$ and the minimum norm $36$ from the lattice ${\bf L}$ and the ternary code ${\bf C}$. The center density of the lattice ${\bf W}_{12}$ is $\frac{1}{48}$. The kissing number $k_1({\bf W}_{12})=3(12+12)+ 8 (3+4+4\cdot2 +(4+1) \cdot 6)=432$ of the lattice ${\bf W}_{12}$ can be calculated similarly as the previous sections.\\

A similar argument as the proof of Theorem 2.2 gives us the following result.\\

{\bf Lemma 8.1.} {\em For each codeword ${\bf c} \in {\bf C}' $ not in ${\bf C}$ and any vector ${\bf V}=({\bf v}_1,{\bf v}_2,{\bf v}_3)$ satisfying  ${\bf v}_1 \equiv {\bf v}_2 \equiv {\bf v}_3 \equiv {\bf c}$ $mod$ $3$, in the $12$ dimensional Euclid space ${\bf R}^{12}$ spanned by the lattice ${\bf W}_{12}$, the distance of the vector ${\bf V}$ to the lattice ${\bf W}_{12}$ is at least $\sqrt{30}$.}\\

Let ${\bf t}=(0,0,\frac{3}{2},\frac{3}{2},-\frac{3}{2},-\frac{3}{2},0,0)$ and ${\bf U}_1=({\bf t},-{\bf t},{\bf 0}),{\bf U}_2=(0,{\bf t},-{\bf t}), {\bf U}_3\\={\bf U}_1+{\bf U}_2=({\bf t},{\bf 0},-{\bf t})$. Here it is clear $2{\bf U}_i \in {\bf W}_{12}$. We have the following result.\\

{\bf Lemma 8.2.} {\em For each codeword ${\bf c} \in {\bf C}' $ not in ${\bf C}$ and a vector ${\bf V}= ({\bf v}_1,{\bf v}_2,{\bf v}_3)$ satisfying
${\bf v}_1 \equiv {\bf v}_2 \equiv {\bf v}_3 \equiv {\bf c}$ $mod$ $3$, in the $12$ dimensional space ${\bf R}^{12}$ spanned by the
lattice ${\bf W}_{12}$, the distance of the vector ${\bf V}$ to ${\bf U}_i$ is at least $\sqrt{30}$.}\\

{\bf Proof.} The Hamming weight of the codeword ${\bf c}$ is $5$ or
$8$ since it is not in the $[8, 2, 6]_3$ ternary code. If the
Hamming weight of ${\bf c}$ is $8$, the Euclid norm of the vector
${\bf V}-{\bf U}_i$ is at least $4 \cdot \frac{3}{2}+4 \cdot 6=30$.
If the Hamming weight of the codeword ${\bf c}$ is $5$ and there are
two elements of the support of ${\bf c}$ in the set $\{3,4,5,6\}$,
the Euclid norm of the vector ${\bf V}-{\bf U}_i$ is at least $2
\cdot \frac{9}{2}+ 2 \cdot \frac{3}{2}+ 3 \cdot 6=30$. If the
Hamming weight of the codeword ${\bf c}$ is $5$ and there are three
elements of the support of ${\bf c}$ in the set $\{3,4,5,6\}$, in
these three row positions of each vector ${\bf v}_i$, there is at
least  one $\pm2$ or $\pm3$ at the row corresponding to the zero row
of ${\bf U}_i$. Thus the Euclid norm of the vector ${\bf V}-{\bf
U}_i$
is at least $2 \cdot 6+ \frac{9}{2}+ 2 \cdot\frac{3}{2}+ (4+\frac{1}{4}+\frac{25}{2})=30$. The conclusion is proved.\\

We construct a new $14$ dimensional lattice ${\bf T}_{14}$ in the $14$ dimensional Euclid space ${\bf R}^{14} \subset {\bf R}^9 \oplus {\bf R}^9 \oplus {\bf R}^9$ with real coordinates $(x_{ij})_{1\leq i \leq 3, 1 \leq j \leq 9}$, where $x_{1j}+x_{2j}+x_{3j}=0$ for each $j=1,...,9$ and $x_{i1}+x_{i2}-x_{i7}-x_{i8}=0$,  $x_{i3}-x_{i4}-x_{i5}+x_{i6}=0$ for $i=1,2,3$. The lattice ${\bf T}_{14}$ is the linear span with integer coefficients of the lattice vectors in ${\bf W}_{12}$ and the following two vectors ${\bf X}$ and ${\bf Y}$.

$$
\left(
\begin{array}{ccccccccccc}
2&0&\frac{1}{2}&\frac{1}{2}&-\frac{3}{2}&-\frac{3}{2}&1&1&\sqrt{3}\\
-1&0&\frac{1}{2}&\frac{1}{2}&\frac{3}{2}&\frac{3}{2}&-2&1&-\sqrt{3}\\
-1&0&-1&-1&0&0&1&-2&0\\
\end{array}
\right)
$$

$$
\left(
\begin{array}{ccccccccccc}
2&0&-1&-1&0&0&1&1&0\\
-1&0&\frac{7}{2}&\frac{7}{2}&-\frac{3}{2}&-\frac{3}{2}&-2&1&\sqrt{3}\\
-1&0&-\frac{5}{2}&-\frac{5}{2}&\frac{3}{2}&\frac{3}{2}&1&-2&-\sqrt{3}\\
\end{array}
\right)
$$

Let ${\bf T}$ be an integer vector of the following form.\\
$$
\left(
\begin{array}{ccccccccccc}
2&0&-1&-1&0&0&1&1\\
-1&0&2&2&0&0&-2&1\\
-1&0&-1&-1&0&0&1&-2\\
\end{array}
\right)
$$
The vector ${\bf T}$ is in the real space spanned by the lattice ${\bf W}_{12}$ and $3{\bf T} \in {\bf W}_{12}$ . The residue class of this vector ${\bf T}$ is a Hamming weight $5$ codeword $(-1,0,-1,-1,0,0,1,1) \in {\bf C}'$. Here it should be noted that the component ${\bf X}_{12345678}$ at $12345678$ positions of the vector ${\bf X}$ is the sum of the vector ${\bf T}$ and the vector ${\bf U}_1$ and the component ${\bf Y}_{(12345678)}$ at $12345678$ positions of the vector ${\bf Y}$ is the sum of the vector ${\bf U}_2$ and the vector ${\bf T}$. We have ${\bf X}_{(12345678)}={\bf T}+{\bf U}_1$
and ${\bf Y}_{(12345678)}={\bf T}+{\bf U}_2$.\\

{\bf Lemma 8.3.} {\em For any integer pair $(m,n) \in {\bf Z}^2$ satisfying $(m,n) \neq (0,0)$, the distance of $m{\bf X}+n{\bf Y}$ to any
lattice vector ${\bf V}$ in ${\bf W}_{12}$ is at least $6$.}\\

{\bf Proof.} First of all the the component of the difference of
$m{\bf X}+n{\bf Y}-{\bf V}$ at the position $9$ is of the form $(m
\sqrt{3},n\sqrt{3}-m\sqrt{3},-n\sqrt{3})$. Its Euclid norm
$6(m^2+n^2-mn)=6((m-\frac{n}{2})^2+\frac{3n^2}{4})$ is smaller than
$36$ if and only if $(m,n)=\pm(1,0)$, $(m,n)=\pm(0,1)$,
$(m,n)=\pm(1,1)$, $(m,n)=\pm(1,-1)$, $(m,n)=\pm(2,1)$,
$(m,n)=\pm(1,2)$ or $(m,n)=\pm(2,2)$. In the last case it can be
checked that the Euclid norm of the difference is at least
$24+30=54$.
In the first three cases the conclusion follows from Lemma 8.1 and Lemma 8.2 directly. In the last three cases the component of the difference vector $m{\bf X}+n{\bf Y}-{\bf V}$ at position $9$ is of the form $(\sqrt{3},-2\sqrt{3},\sqrt{3})$ (or equivalent other forms) and its Euclid norm is $18$. The component of the difference vector $m{\bf X}+n{\bf Y}-{\bf V}$ at the positions $(12345678)$ is ${\bf U}_1-{\bf U}_2=({\bf t},-2{\bf t},{\bf t})$. It is easy to check that the Euclid norm of the difference of $({\bf t},-2{\bf t},{\bf t})$ to any lattice vector in the lattice ${\bf W}_{12}$ is at least $18$. The conclusion is proved.\\

The volume of the lattice ${\bf T}_{14}$ is $3^{13} \cdot 16 \cdot \sqrt{6} \cdot \sqrt{\frac{9}{2}}=3^{14} \cdot 16\sqrt{3}$. Thus the center density of ${\bf T}_{14}$ is $\frac{1}{16 \sqrt{3}}$. \\

The Gram matrix of the lattice ${\bf T}_{14}$ (with repect to one base) is as follows.\\

$$
\left(
\begin{array}{cccccccccccccccc}
36&18&18&9&0& 0&0&0&0&0&-18&0&9&9\\
18&36&9&18&0&0&0&0&0&0&-18&0&9&9\\
18&9&36&18&9 &0&0&0&0&0&-9&-9&18&18\\
9&18&18&36&18&0&0&0&0&0&-9& -9 &9& 9\\
0&0&9&18 &36 &0&0&0&0 &0&0&-18&-9&-9\\
0&0&0 &0&0&36&0&0&9&18&-18&0&9&9\\
0&0&0&0&0&0&36&18&-18&-9 &18&0&-18&9\\
0&0&0&0&0&0&18&36&-9&-18&18&0&-9&-9\\
0&0&0&0&0&9 &-18 &-9&36&18 &-18&18&9&-18\\
0&0&0&0&0&18&-9&-18&18&36&-18&18&9&9\\
-18&-18&-9&-9&0&-18&18&18&-18&-18&36&0&-18&0\\
0&0&-9&-9&-18&0&0&0&18&18&0&36&0&0\\
9&9&18&9&-9&9&-18&-9& 9&9&-18&0&36&18\\
9&9&18&9&-9&9&9&-9&-18&9&0&0&18&72\\
\end{array}
\right)
$$

{\bf Corollary 8.4} {\em There is a $13$ dimensional lattice with the center density $\frac{1}{16\sqrt{6}}$ and the kissing number $726$. There is a family of $14$ dimenisonal lattices ${\bf T}_{14,t}$ with the kissing number $k_1({\bf T}_{14,t})=432$ and the center density tending to $\frac{1}{16 \sqrt{3}}$.}\\

{\bf Proof.} The $13$ dimensional lattice can be generated by ${\bf W}_{12}$ and the vector ${\bf X}$. If  we adjust $\sqrt{3}$ in the ${\bf X}$ and ${\bf Y}$ to a little bigger real number and set them tending to $\sqrt{3}$ we get the family of the $14$ dimensional lattices.\\

From the above Gram matrix we can get an isometry to the lattice ${\bf \Lambda}_{14}$ by Magma.\\

The following is a simpler recovery of the $14$ dimensional lattice ${\bf \Lambda}_{14}$.\\

Let ${\bf C}$ be the following $[8, 3, 3]_3$ ternary code with the following generator matrix.\\

$$
\left(
\begin{array}{cccccccccccccccc}
1&1&1&0&1&1&1&0\\
0&1&-1&1&0&1&-1&1\\
1&1&1&0&0&0&0&0\\
\end{array}
\right)
$$

There are $4$ Hamming weight $3$ codewords $\pm(1,1,1,0,0,0,0,0),\pm(0,0,0,\\0,1,1,1,0)$ and $22$ Hamming weight $6$ codewords in ${\bf C}$. Among these 22 codewors of Hamming weight 6, 10 of them are symmetric, these are $\pm(1,1,1,0,\\-1,-1,1,0)$ and the eight non-zero codewords of the $[8,2,6]_3$ code generated by the first two rows, the other twelve codwords of Hamming weight six are not symmetric because  the 2rd rows are invovled. Here we should note that $(0,0,0,1,0,0,0,-1)$ is in the dual code of ${\bf C}$.\\

Let $L \subset ({\bf Z}^8 \oplus {\bf Z}^8 \oplus {\bf Z}^8)$ be a $14$ dimensional lattice defined by the column conditions $x_{1j}+x_{2j}+x_{3j}=0$, where $j=1,...,8$ and $x_{ij}, i=1,2,3,j=1,...,8$'s are the coordinates of the lattice ${\bf Z}^8 \oplus {\bf Z}^8 \oplus {\bf Z}^8$ and $x_{i4}=x_{i8}$ for $i=1,2,3$. This is a lattice with the volume $3^{3.5} \cdot 2$.
From Theorem 2.3 we get a $14$ dimensional lattice ${\bf S'}_{14}$ with the minimum norm $18$ and the volume $3^{14.5} \cdot 2$. Its center density is $\frac{1}{256\sqrt{3}}$. \\

The following vectors ${\bf X}_1,...,{\bf X}_7$ \\

$$
\left(
\begin{array}{cccccccccccccccc}
\frac{3}{2}&\frac{3}{2}&0&0&\frac{3}{2}&\frac{3}{2}&0&0\\
-\frac{3}{2}&-\frac{3}{2}&0&0&-\frac{3}{2}&-\frac{3}{2}&0&0\\
0&0&0&0&0&0&0&0\\
\end{array}
\right)
$$

$$
\left(
\begin{array}{cccccccccccccccc}
0&\frac{3}{2}&\frac{3}{2}&0&0&\frac{3}{2}&\frac{3}{2}&0\\
0&-\frac{3}{2}&-\frac{3}{2}&0&0&-\frac{3}{2}&-\frac{3}{2}&0\\
0&0&0&0&0&0&0&0\\
\end{array}
\right)
$$

$$
\left(
\begin{array}{cccccccccccccccc}
\frac{3}{2}&0&\frac{3}{2}&0&\frac{3}{2}&0&\frac{3}{2}&0\\
-\frac{3}{2}&0&-\frac{3}{2}&0&-\frac{3}{2}&0&-\frac{3}{2}&0\\
0&0&0&0&0&0&0&0\\
\end{array}
\right)
$$

$$
\left(
\begin{array}{cccccccccccccccc}
\frac{3}{2}&0&0&\frac{3}{2}&\frac{3}{2}&0&0&\frac{3}{2}\\
-\frac{3}{2}&0&0&-\frac{3}{2}&-\frac{3}{2}&0&0&-\frac{3}{2}\\
0&0&0&0&0&0&0&0\\
\end{array}
\right)
$$

$$
\left(
\begin{array}{cccccccccccccccc}
0&\frac{3}{2}&0&\frac{3}{2}&0&\frac{3}{2}&0&\frac{3}{2}\\
0&-\frac{3}{2}&0&-\frac{3}{2}&0&-\frac{3}{2}&0&-\frac{3}{2}\\
0&0&0&0&0&0&0&0\\
\end{array}
\right)
$$

$$
\left(
\begin{array}{cccccccccccccccc}
0&0&\frac{3}{2}&\frac{3}{2}&0&0&\frac{3}{2}&\frac{3}{2}\\
0&0&-\frac{3}{2}&-\frac{3}{2}&0&0&-\frac{3}{2}&-\frac{3}{2}\\
0&0&0&0&0&0&0&0\\
\end{array}
\right)
$$

$$
\left(
\begin{array}{cccccccccccccccc}
\frac{3}{2}&\frac{3}{2}&\frac{3}{2}&\frac{3}{2}&\frac{3}{2}&\frac{3}{2}&\frac{3}{2}&\frac{3}{2}\\
-\frac{3}{2}&-\frac{3}{2}&-\frac{3}{2}&-\frac{3}{2}&-\frac{3}{2}&-\frac{3}{2}&-\frac{3}{2}&-\frac{3}{2}\\
0&0&0&0&0&0&0&0\\
\end{array}
\right)
$$
and the following vector ${\bf Y}$ \\

$$
\left(
\begin{array}{cccccccccccccccc}
0&0&0&0&0&0&0&0\\
-\frac{3}{2}&-\frac{3}{2}&-\frac{3}{2}&-\frac{3}{2}&-\frac{3}{2}&-\frac{3}{2}&-\frac{3}{2}&-\frac{3}{2}\\
\frac{3}{2}&\frac{3}{2}&\frac{3}{2}&\frac{3}{2}&\frac{3}{2}&\frac{3}{2}&\frac{3}{2}&\frac{3}{2}\\
\end{array}
\right)
$$
are needed. It is easy to check that the above $8$ vectors are in the real space ${\bf R}^{14}$ spanned by the lattice ${\bf S'}_{14}$.\\

We consider the $16$ vectors of the form ${\bf 0},{\bf X}_1,{\bf X}_2,{\bf X}_3$,${\bf X}_4,{\bf X}_5,{\bf X}_6,{\bf X}_7$,${\bf Y}, {\bf Y}-{\bf X}_1,{\bf Y}-{\bf X}_2,{\bf Y}-{\bf X}_3$,${\bf Y}-{\bf X}_4,{\bf Y}-{\bf X}_5,{\bf Y}-{\bf X}_6,{\bf Y}-{\bf X}_7$.  These vectors are in the real space sapnned by the lattice ${\bf S'}_{14}$. The union of the sixteen translates of the above lattice ${\bf S'}_{14}$ is a lattice ${\bf S}_{14}$ with the volume $\frac{3^{14.5}}{8}$. Here it should be noted that $2{\bf X}_i, 2{\bf Y} \in {\bf L} \subset {\bf S'}_{14}$.\\

{\bf Theorem 8.5.} {\em The minimum norm of the lattice ${\bf S}_{14}$ is $18$.}\\

{\bf Proof.}  Similar to the proof of Theorem 51.\\

Thus we get a $14$ dimensional lattice ${\bf S}_{14}$ with the center density $\frac{(\frac{18}{4})^7}{\frac{3^{14.5}}{8}}=\frac{1}{16\sqrt{3}}$. \\

By using a base we can get an isometry to the laminated lattice ${\bf \Lambda}_{14}$ from Magma.\\

\section{Some other easy lattices}

From a similar construction of ${\bf T}_{10}$ and ${\bf T}_{12}$, by using the ternary $[10,4,6]_3$ code we can get a $18$ dimensional lattice
${\bf T}_{18,1}$ with the center density $\frac{1}{10\sqrt{3}}$ and the kissing number
$3 \cdot 10 \cdot9 +18 \cdot 6 \cdot 6 +42 (3+3 \cdot2 \cdot 3 \cdot 3+3\cdot 3 \cdot 2\cdot2 )=5796$.
This lattice has the same center density as the  "physical" lattice $dim18kiss5820min4det300$ in \cite{Nebe} found
by Y. Kallus (\cite{Kallus}) but with a smaller kissing number. From the $[11, 5, 6]_3$ ternary code
(the weight distribution $A_0=1,A_6=132,A_9=110$) and two disjoint codewords of Hamming weight $5$ and $6$ in its dual,
a $18$ dimensional lattice ${\bf T}_{18,2}$ with the center density $\frac{1}{10 \sqrt{3}}$ and the kissing number
$3(5\cdot 4+6 \cdot5)+42 \cdot 36 +90 \cdot 45=150+1512+4050=5712$ can be constructed.
This lattice has the same center density as the lattice $dim18kiss5820min4det300$ in \cite{Nebe}
found by Y. Kallus (\cite{Kallus}),  but with a smaller kissing number.  From the ternary $[11, 5, 6]_3$ ternary code
a $20$ dimensional lattice ${\bf T}_{20}$ with the center density $\frac{1}{11}$ can also be constructed.
These lattices are worse than the Kappa lattices ${\bf K}_{18}$ and ${\bf K}_{20}$. \\

From the $[9, 3, 6]_3$ ternary code in section 4 we can get a $14$ dimensional lattice with the center density $\frac{1}{18\sqrt{3}}$ and the kissing number $1248$. That is the Kappa lattice ${\bf K}_{14,2}$ (\cite{Plesken,Nebe}). Another lattice ${\bf K}_{14}$ with the center density $\frac{1}{18\sqrt{3}}$ and the kissing number $1242$ (\cite{Plesken,Nebe}) can also be constructed similarly. Actually from any two disjoint codwords with Hamming weight three and Hamming weight six in the dual code such a $14$ dimensional lattice can be constructed. It seems they are all isomorphic to the Kappa lattices of dimension $14$. If we use three codewords in the dual of the ternary $[10, 4, 6]_3$ code with Hamming weight $4$, $6$ and $4$ a $14$ dimensional lattice ${\bf U}_{14}$ with the center density $\frac{1}{18\sqrt{3}}$ can be constructed. This lattice would have a "smaller" kissing number. We guess that it is not one of the Kappa lattices of dimension $14$. If we use three codewords in the dual of the ternary $[10, 4, 6]_3$ code with Hamming weight $5$, $5$ and $6$ a $14$ dimensional lattice ${\bf W}_{14}$ with the center density $\frac{1}{20\sqrt{3}}$ can be constructed. It is worser than the laminated lattices and the Kappa lattices of dimension $14$, however it is certainly better than the "physical" lattice $dim14kiiss774min8det25021632$ found by Y. Kallus (\cite{Kallus}) in \cite{Nebe}.\\

From the $[8,2,6]_3$ ternary code used in section 3 we can get a $12$ dimensional lattice ${\bf S}_{12}$ with the center density $\frac{1}{36}$ and the kissing number $456$ can be constructed. This lattice is worser than the Laminated lattices of dimension $12$ and the Coxeter-Todd lattice, but better than the " physical" lattice $dim12kiss462min8det5619712$ in (\cite{Kallus},\cite{Nebe}). By eqauting two column coordinates we have a $11$ dimensional lattice ${\bf S}_{11}$ with the center density $\frac{1}{36}$.\\

From the $[10, 4,6]_3$ ternary code, we can construct a $16$ dimensional lattices with the center density $\frac{1}{24}$ from any two disjoint weight four and weight six codewords in the dual code. This lattice has the same center density as the Kappa lattices of dimension $16$. One lattice ${\bf S}_{16}$ of them has the kissing number $3(4\cdot3+6 \cdot 5)+18 \cdot 6 \cdot 6+42 \cdot 3 \cdot 15=126+360+2250=2664$. This is not the Kappa lattice ${\bf K}_{16}$ of dimension $16$ with center density $\frac{1}{24}$ and the kissing number $2772$ (\cite{Plesken,Nebe}). By equating two column coordinates of this lattice we have a $15$ dimensional lattice ${\bf S}_{15}$ with the center density $\frac{1}{12\sqrt{6}}$. It is as dense as the  Kappa lattices ${\bf K}_{15}$ in \cite{Plesken,Nebe}. However this lattice is not ${\bf K}_{15}$ since it is a sub-lattice ${\bf S}_{16}$. If we use two disjoint codewords of Hamming weight $5$ in the dual of the ternary $[10,4,6]_3$ code we get a $16$ dimensional lattice with the center density $\frac{1}{25}$ and the kissing number $2712$. The density of this lattice is slightly smaller than the "physical" lattice $dim16kiss2400min6det390625$ (\cite{Kallus,Nebe}).\\

By using the ternary $[11, 5, 6]$ linear code we can get a $22$ dimensional lattice with the center density $\frac{1}{4\sqrt{3}}$ and the kissing number $11 \cdot 10 \cdot 4 \cdot 2+ 132(3+15 \cdot 6 +15 \cdot 6)=24816$ (as the lattice ${\bf T'}_{24}$ in section 6).\\

The following is the base of a $10$ dimensional lattice ${\bf H}_{10}$ with the center density $\frac{1}{32}$ and the kissing number $276$.  This lattice is the $10$ dimensional lattice from the binary linear $[10, 5, 4]_2$ code by using the "Construction A". This seems indicating that the ternary construction of lattice proposed in this paper is quite general.\\
$v1 := V![-1, -1, -1, -1, -1, -1, 3, 3,   \frac{1}{2}, \frac{1}{2}, \frac{1}{2}, \frac{1}{2}, \frac{1}{2},  \frac{1}{2}, -\frac{3}{2},  -\frac{3}{2},   \frac{1}{2}, \frac{1}{2}, \frac{1}{2},  \frac{1}{2},\frac{1}{2}, \frac{1}{2}, \\ -\frac{3}{2}, -\frac{3}{2}]$\\
$v2 := V![\frac{1}{2}, \frac{1}{2}, \frac{1}{2},   \frac{1}{2}, \frac{1}{2}, \frac{1}{2}, -\frac{3}{2}, -\frac{3}{2},   -1, -1, -1, -1, -1, -1, 3, 3,   \frac{1}{2},  \frac{1}{2}, \frac{1}{2},   \frac{1}{2}, \frac{1}{2},  \frac{1}{2}, \\ -\frac{3}{2},  -\frac{3}{2}]$\\
$v3 := V![1, -1, -1,0,0,1,1,-1,   -2,2,2,0, 0,-2, -2,2,      1,-1, -1,0, 0,1,1,-1]$\\
$v4 :=V![0,0,1,-1,1,-1,1,-1,  0,0,-2,2,-2,2,-2,2,   0,0,1,-1,1,-1,1,-1]$\\
$v5 :=V![3,3,-3,0,-3,0,0,0,    -3,-3,3,0,3,0,0,0, 0,0,0,0,0,0,0,0]$\\
$v6 :=V![0,0,0,0, 0,0,0,0,     3,3,-3,0,-3,0,0,0,  -3,-3,3,0,3,0,0,0]$\\
$v7 := V![3,   -3,   0,   0,   0,   0,  0,   0,      -3,  3,   0,   0,    0,   0,   0,   0,    0,   0,   0,   0,   0,  0,   0,   0]$\\
$v8 := V![0,     0,   0,   0,   0,   0,  0,   0,    3 , -3, -3,   0,    0,   3,   0,   0,     -3,   3,   3,   0,   0,   -3,  0,   0]$\\
$v9 := V![0,   0,   0,   0,   0,   0,  0,   0,      -3,  -3,   0,   0,    0,   0,   3,   3,    3,  3,    0,   0,   0,   0,  -3,   -3]$\\
$v10 :=V![3,    3,    0,  0,   0,   0,  -3,   -3,     -3,   -3,   0,   0,   0,   0,    3,   3,   0,   0,   0,   0,   0,  0,  0,  0]$\\

Its Gram matrix is as follows.\\
$$
\left(
\begin{array}{ccccccccccccccccc}
36& -18 &  0 &  0  & 0 &  0 &  0 &  0 &  0 &-36\\
-18&  36 &  0&   0 &  0&   0&   0&   0 & 36&  36\\
 0&   0&  36&   0&   9&  -9&  18& -36&   0&   0\\
0 &  0&   0 & 36& -18&  18 &  0 & 18 &  0&   0\\
0 &  0&   9& -18&  72& -36 &  0&  -9&  18&  36\\
0 &  0&  -9&  18& -36&  72&   0&  18& -36& -18\\
0&   0&  18 &  0&   0&   0&  36& -18&   0&   0\\
0&   0& -36&  18 & -9&  18& -18 & 72&   0&   0\\
0&  36&   0&   0&  18& -36&   0 &  0 & 72 & 36\\
-36&  36&   0&   0&  36& -18&   0&   0&  36&  72\\
\end{array}
\right)
$$

The theta series of the lattice ${\bf H}_{10}$ is as follows.\\

$1 + 276q^4 + 768q^6 + 4020q^8 + 6144q^{10} + 20416q^{12} + 23040q^{14} +\\
    65844q^{16} + 61440q^{18} + 160488q^{20}+ 140544q^{22} + 327616q^{24} +\\
    276480q^{26} + 612480q^{28} + 480768q^{30} + 1047348q^{32} + 798720q^{34} +
    1665876q^{36} + 1251072q^{38} + 2565672q^{40} + ...$.\\

The lattice generated by the following vectors is $14$ dimensional lattice with the center density $\frac{1}{18\sqrt{3}}$ and the kissing number $1248$. By Magma we can check that this lattice is the lattice ${\bf K}_{14,2}$ (\cite{Nebe,Plesken}).\\
$v1 := V![-1,-1,-1,-1,-1,-1,3,3, \frac{1}{2},\frac{1}{2},\frac{1}{2},\frac{1}{2},\frac{1}{2},\frac{1}{2}, -\frac{3}{2},-\frac{3}{2}, \frac{1}{2},\frac{1}{2},\frac{1}{2}, \frac{1}{2},  \frac{1}{2}, \frac{1}{2}, \\ -\frac{3}{2},-\frac{3}{2}]$\\
$v2 := V![\frac{1}{2},\frac{1}{2},\frac{1}{2},\frac{1}{2},\frac{1}{2},\frac{1}{2}, -\frac{3}{2},-\frac{3}{2}, -1,-1,-1,-1, -1,-1,3,3, \frac{1}{2},\frac{1}{2},\frac{1}{2},  \frac{1}{2}, \frac{1}{2}, \frac{1}{2},\\  -\frac{3}{2},-\frac{3}{2}]$\\
$v3 :=V![3,0,-3,0,0,0,0,0, -3,0,3,0,0,0,0,0, 0,0,0,0,0,0,0,0]$\\
$v4 :=V![0,0,0,0, 0,0,0, 0, 3,0,-3,0,0,0,0,0,  -3,0,3,0,0,0,0,0]$\\
$v5 := V![3,0,0,-3,0,0, 0,0,  -3,0,0,3,0,0, 0,0,    0,0,0,0,   0,0,0,0]$\\
$v6 := V![0,0,0,0, 0,0,0,0, 3,0,0,-3,0,0, 0,0,  -3,0,0,3,0,0, 0,0]$\\
$v7 :=V![3,0,0,0,-3,0,0,0,    -3,0,0,0,3,0,0,0, 0,0,0,0,0,0,0,0]$\\
$v8 :=V![0,0,0,0, 0,0,0,0,    3,0,0,0,-3,0,0,0,  -3,0,0,0,3,0,0,0]$\\
$v9 := V![3, 0, 0, 0, 0, -3,0, 0,      -3, 0, 0, 0,0,   3,   0,   0,     0,0,0,0,   0,0,0,0]$\\
$v10 := V![0,0,0,0, 0,0,0,0, 3, 0, 0, 0, 0, -3,0, 0,      -3, 0, 0, 0,0,   3,   0,   0]$\\
$v11 := V![3,0,0,0,   0,0,-3,0,      -3,0,0,0,0,0,   3,0,    0,  0,    0,   0,   0,   0,  0,0]$\\
$v12 := V![0,0,0,0, 0,0,0,0,   3,0,0,0,   0,0,-3,0,      -3,0,0,0, 0,0,   3,0]$\\
$v13 :=V![3,0,    0,0,0,0,  0,-3,     -3,0,   0,0,0,0,0,3,   0,   0,   0,   0,   0,  0,  0,  0]$\\
$v14 :=V![0,0,0,0, 0,0,0,0, 3,0,    0,0,0,0,  0,-3,     -3,0,   0,0,0,0,0,3]$\\

Its theta series is $1 + 1248q^4 + 11808q^6 + 73062q^8 + 263520q^{10} + 811032q^{12} + 2019168q^{14}
    + 4542864q^{16} + 8999136q^{18} + 17694288q^{20} + \cdots$\\

The lattice generated by the following vectors  is a $13$ dimensional sub-lattice of the above lattice ${\bf K}_{14,2}$. Its center density is $\frac{1}{36}$ and the kissing number is $780$. This lattice can be compared with the lattice "dim13kis720min8det\\10616832" in \cite{Kallus,Nebe}.\\
$v1 := V![-1,-1,-1,-1,-1,-1,3,3, \frac{1}{2},\frac{1}{2},\frac{1}{2},\frac{1}{2},\frac{1}{2},\frac{1}{2}, -\frac{3}{2},-\frac{3}{2}, \frac{1}{2},\frac{1}{2},\frac{1}{2}, \frac{1}{2},  \frac{1}{2}, \frac{1}{2}, \\ -\frac{3}{2},-\frac{3}{2}]$\\
$v2 := V![\frac{1}{2},\frac{1}{2},\frac{1}{2},\frac{1}{2},\frac{1}{2},\frac{1}{2}, -\frac{3}{2},-\frac{3}{2}, -1,-1,-1,-1, -1,-1,3,3, \frac{1}{2},\frac{1}{2},\frac{1}{2},  \frac{1}{2}, \frac{1}{2}, \frac{1}{2},\\  -\frac{3}{2},-\frac{3}{2}]$\\
$v3 :=V![0,0,0,0, 0,0,0, 0, 3,0,-3,0,0,0,0,0,  -3,0,3,0,0,0,0,0]$\\
$v4 := V![3,0,0,-3,0,0, 0,0,  -3,0,0,3,0,0, 0,0,    0,0,0,0,   0,0,0,0]$\\
$v5 := V![0,0,0,0, 0,0,0,0, 3,0,0,-3,0,0, 0,0,  -3,0,0,3,0,0, 0,0]$\\
$v6 :=V![3,0,0,0,-3,0,0,0,    -3,0,0,0,3,0,0,0, 0,0,0,0,0,0,0,0]$\\
$v7 :=V![0,0,0,0, 0,0,0,0,    3,0,0,0,-3,0,0,0,  -3,0,0,0,3,0,0,0]$\\
$v8 := V![3, 0, 0, 0, 0, -3,0, 0,      -3, 0, 0, 0,0,   3,   0,   0,     0,0,0,0,   0,0,0,0]$\\
$v9 := V![0,0,0,0, 0,0,0,0, 3, 0, 0, 0, 0, -3,0, 0,      -3, 0, 0, 0,0,   3,   0,   0]$\\
$v10 := V![3,0,0,0,   0,0,-3,0,      -3,0,0,0,0,0,   3,0,    0,  0,    0,   0,   0,   0,  0,0]$\\
$v11 := V![0,0,0,0, 0,0,0,0,   3,0,0,0,   0,0,-3,0,      -3,0,0,0, 0,0,   3,0]$\\
$v12 :=V![3,0,    0,0,0,0,  0,-3,     -3,0,   0,0,0,0,0,3,   0,   0,   0,   0,   0,  0,  0,  0]$\\
$v13 :=V![0,0,0,0, 0,0,0,0, 3,0,    0,0,0,0,  0,-3,     -3,0,   0,0,0,0,0,3]$\\

Its theta series is $1 + 780q^{36} + 5784q^{54} + 32022q^{72} + 100264q^{90} + \cdots$.\\

Actually the above 14 dimensional lattice ${\bf K}_{14,2}$ has many $13$ dimensional sub-lattices with the same parameters.\\

The above lattice ${\bf K}_{14,2}$ also has the following $12$ dimensional sub-lattice. Its center density is $\frac{1}{36}$ and the kissing number is $552$. \\
$v1 := V![-1,-1,-1,-1,-1,-1,3,3, \frac{1}{2},\frac{1}{2},\frac{1}{2},\frac{1}{2},\frac{1}{2},\frac{1}{2}, -\frac{3}{2},-\frac{3}{2}, \frac{1}{2},\frac{1}{2},\frac{1}{2}, \frac{1}{2},  \frac{1}{2}, \frac{1}{2}, \\ -\frac{3}{2},-\frac{3}{2}]$\\
$v2 := V![\frac{1}{2},\frac{1}{2},\frac{1}{2},\frac{1}{2},\frac{1}{2},\frac{1}{2}, -\frac{3}{2},-\frac{3}{2}, -1,-1,-1,-1, -1,-1,3,3, \frac{1}{2},\frac{1}{2},\frac{1}{2},  \frac{1}{2}, \frac{1}{2}, \frac{1}{2},\\  -\frac{3}{2},-\frac{3}{2}]$\\
$v3 :=V![0,0,0,0, 0,0,0, 0, 3,0,-3,0,0,0,0,0,  -3,0,3,0,0,0,0,0]$\\
$v4 := V![0,0,0,0, 0,0,0,0, 3,0,0,-3,0,0, 0,0,  -3,0,0,3,0,0, 0,0]$\\
$v5 :=V![3,0,0,0,-3,0,0,0,    -3,0,0,0,3,0,0,0, 0,0,0,0,0,0,0,0]$\\
$v6 :=V![0,0,0,0, 0,0,0,0,    3,0,0,0,-3,0,0,0,  -3,0,0,0,3,0,0,0]$\\
$v7 := V![3, 0, 0, 0, 0, -3,0, 0,      -3, 0, 0, 0,0,   3,   0,   0,     0,0,0,0,   0,0,0,0]$\\
$v8 := V![0,0,0,0, 0,0,0,0, 3, 0, 0, 0, 0, -3,0, 0,      -3, 0, 0, 0,0,   3,   0,   0]$\\
$v9 := V![3,0,0,0,   0,0,-3,0,      -3,0,0,0,0,0,   3,0,    0,  0,    0,   0,   0,   0,  0,0]$\\
$v10 := V![0,0,0,0, 0,0,0,0,   3,0,0,0,   0,0,-3,0,      -3,0,0,0, 0,0,   3,0]$\\
$v11 :=V![3,0,    0,0,0,0,  0,-3,     -3,0,   0,0,0,0,0,3,   0,   0,   0,   0,   0,  0,  0,  0]$\\
$v12 :=V![0,0,0,0, 0,0,0,0, 3,0,    0,0,0,0,  0,-3,     -3,0,   0,0,0,0,0,3]$\\

Its theta series is $1 + 552q^{36} + 3048q^{54} + 15858q^{72} + 41544q^{90} + 116436q^{108} + 228312q^{126}
    + 491376q^{144} + 787200q^{162} + 1527300q^{180} + 2122128q^{198} + \cdots$.\\

The following is a $10$ dimensional lattice with the center density $\frac{3\sqrt{3}}{160}=\frac{1}{17.77777...\sqrt{3}}=\frac{1}{30.7920143565...}=0.03247....$. The kissing number is $272$. It is smaller than the kissing number $276$ of the Kappa lattice ${\bf K}_{10}$ and bigger than the kissing number $270$ of the Martinet lattice ${\bf K}_{10}'$ (\cite{Mar}, page 290). Though the center density is slightly bigger than the center density $\frac{1}{18\sqrt{3}}$ of these two lattices.\\

Let ${\bf C}$ be the ternary Hamming $[4, 2,3]_3$ linear code. Let $(x)_{1 \leq i \leq 3, 1 \leq j \leq 18}$ be the coordinates of the $54$ dimensional real space ${\bf R}^{18} \oplus {\bf R}^{18} \oplus {\bf R}^{18}$. The lattice ${\bf L}_{10}' \subset {\bf Z}^{18} \oplus {\bf Z}^{18} \oplus {\bf Z}^{18}$ is defined by the following conditions.\\
1) $x_{1j}+x_{2j}+x_{3j}=0$, $j=1,...,18$;\\
2) $x_{i1}=x_{i2}$, $x_{i3}=x_{i4}$, $x_{i5}=x_{i6}$, $x_{i7}=x_{i8}$, $x_{i9}=x_{i,10}=\cdots=x_{i18}$; $i=1,2,3$;\\
3) $(x_{11}-x_{21}, x_{13}-x_{23},x_{15}-x_{25},x_{17}-x_{27}) \in {\bf C}$.\\

{\bf Proposition 10.1} {\em The volume of this lattice ${\bf L}_{10}'$ is $160 \cdot 3^{4.5}$ and the minimum norm of ${\bf L}_{10}'$ is $12$.}\\

{\bf Proof.}  If we use the parity check matrix for the third condition, it is clear the volume of the lattice ${\bf L}_{10}'$ is $3^{2.5} \cdot  3^{4-2} \cdot 2 ^4 \cdot 10=160 \cdot 3^{4.5}$. \\

For any lattice vector ${\bf V}=(v_{ij})_{1 \leq i \leq 3, 1 \leq j \leq n}, \in {\bf L}_{2n}$, if $(v_{11}-v_{21}, ...., v_{1,18}-v_{2,18})$, $mod$ $3$ is a zero in ${\bf F}_3^{18}$, then the norm of each non-zero column of ${\bf V}$ is at least $6$, and when the norm of a non-zero column is $6$, this column is $\pm(1,1,-2)^{\tau}$ or their permuatations. If the norm of a non-zero column of ${\bf V}$ is not of this form, the norm is at least $24$. In the case that the non-zero columns of ${\bf V}$ are of the form $\pm(1,1,-2)^{\tau}$ or their permuatations, from the second condition we have at least two such columns. Thus the norm of ${\bf V}$ is at least $12$. \\

If $(v_{11}-v_{21}, ...., v_{1,18}-v_{2,18})$, $mod$ $3$ is a non-zero codeword in ${\bf F}_3^{18}$, If the codeword $(v_{11}-v_{21}, v_{13}-v_{23},v_{15}-v_{25},v_{17}-v_{27})$,  $mod$ $3$ is zero codeword in ${\bf C}$, $(v_{19}-v_{29}, v_{1,10}-v_{2,10},..., v_{1,18}-v_{2,18})$, $mod$ $3$ is non-zero in ${\bf F}_3^{10}$. The norm of the part of the last ten columns is at least $2 \cdot10=20$. Otherwise we have at least $2 \cdot 3$ non-zero columns in the part of the first eight columns, the norm of this part is at least $2 \cdot 6=12$.\\

Let ${\bf X}$ be the following vector in the real space spanned by the lattice ${\bf L}_{10}'$. It is clear $3{\bf X} \in {\bf L}_{10}'$.\\

$$
\left(
\begin{array}{ccccccccccccccccccccccccccccccc}
\frac{1}{3}&\frac{1}{3}&\frac{1}{3}&\frac{1}{3}&\frac{1}{3}&\frac{1}{3}&\frac{1}{3}&\frac{1}{3}&\frac{1}{3}&\frac{1}{3}&\frac{1}{3}&\frac{1}{3}&\frac{1}{3}&\frac{1}{3}&\frac{1}{3}&\frac{1}{3}&\frac{1}{3}&\frac{1}{3}\\
\frac{1}{3}&\frac{1}{3}&\frac{1}{3}&\frac{1}{3}&\frac{1}{3}&\frac{1}{3}&\frac{1}{3}&\frac{1}{3}&\frac{1}{3}&\frac{1}{3}&\frac{1}{3}&\frac{1}{3}&\frac{1}{3}&\frac{1}{3}&\frac{1}{3}&\frac{1}{3}&\frac{1}{3}&\frac{1}{3}\\
\frac{-2}{3}&\frac{-2}{3}&\frac{-2}{3}&\frac{-2}{3}&\frac{-2}{3}&\frac{-2}{3}&\frac{-2}{3}&\frac{-2}{3}&\frac{-2}{3}&\frac{-2}{3}&\frac{-2}{3}&\frac{-2}{3}&\frac{-2}{3}&\frac{-2}{3}&\frac{-2}{3}&\frac{-2}{3}&\frac{-2}{3}&\frac{-2}{3}\\
\end{array}
\right)
$$

The union of three translates ${\bf L}_{10}'$, ${\bf X}+{\bf L}_{10}'$ and $-{\bf X}+{\bf L}_{10}'$ is a $10$ dimensional lattice ${\bf L}_{10}$ with the volume $160 \cdot 3^{3.5}$. It is clear the norm of the difference ${\bf X}-{\bf V}$, ${\bf V} \in {\bf L}_{10}'$, is at least $\frac{2}{3} \cdot 18=12$. The center density of the $10$ dimensional lattice ${\bf L}_{10}$ is $\frac{3 \sqrt{3}}{160}$. The kissing number (the number of norm $12$ lattice vectors) in ${\bf L}_{10}$ is $8 \cdot 27 +6 \cdot 4 + 2 \cdot (2  \cdot 2 \cdot 4)=272$. \\

\section{The recovery of the Kappa lattices ${\bf K}_{12}$ (again), ${\bf K}_{18}$, ${\bf K}_{20}$, and a $22$ dimensional lattice ${\bf K}_{22}'$ with the center density $\frac{1}{3\sqrt{3}}$}

Let ${\bf C}$ be the linear ternary $[n, k ,6]_3$ code satisfying the condition that the sum of all coordinates of each codeword is zero. Let $(x_{ij})_{1\leq i, \leq 3, 1 \leq j \leq n}$ be the coordinates of the $3n$ dimensional real space ${\bf R}^n \oplus {\bf R}^n \oplus {\bf R}^n$. We consider the following $2n$ dimensional lattice ${\bf L}_{2n} \subset {\bf Z}^n \oplus {\bf Z}^n \oplus {\bf Z}^n$ defined by the following conditions.\\
1) $x_{1j}+x_{2j}+x_{3j}=0$; \\
2) $x_{11}+x_{12}+\cdots +x_{1n} \equiv 0$, $mod$ $3$;\\
3) $(x_{11}-x_{21}, x_{12}-x_{22}, ..., x_{1n} -x_{2n}) \in {\bf C}$.\\

{\bf Theorem 11.1.} {\em The volume of this $2n$ dimensional lattice ${\bf L}_{2n}$ is $3^{3n/2-k+1}$ and the minimum norm of this lattice is $12$. The kissing number is ${n \choose 2} \cdot 18 +243 A_6$, where $A_6$ is the number of Hamming weight $6$ codewords in ${\bf C}$.}\\

{\bf Proof.} If we use the parity check matrix for the third condition, it is clear the volume of the lattice ${\bf L}_{2n}$ is $3^{n/2} \cdot 3 \cdot 3^{n-k}$. \\

Type 1) For any lattice vector ${\bf V}=(v_{ij})_{1 \leq i \leq 3, 1 \leq j \leq n}, \in {\bf L}_{2n}$, if $(v_{11}-v_{21}, ...., v_{1n}-v_{2n})$, $mod$ $3$ is a zero codeword in ${\bf C}$, then the norm of each non-zero column of ${\bf V}$ is at least $6$, and when the norm of a non-zero column is $6$, this column is $\pm(1,1,-2)^{\tau}$ or their permuatations. If the norm of a non-zero column of ${\bf V}$ is not of this form, the norm is at least $24$. In the case that the non-zero columns of ${\bf V}$ are of the form $\pm(1,1,-2)^{\tau}$ or their permuatations, from the second condition we have at least two such columns. Thus the norm of ${\bf V}$ is at least $12$. \\

Type 2) If $(v_{11}-v_{21}, ...., v_{1n}-v_{2n})$, $mod$ $3$ is a non-zero codeword in ${\bf C}$, there are at least six non-zero columns in ${\bf V}$, since the norm of  each non-zero column is at least $2$, the norm of this lattice vector ${\bf L}_{2n}$ is at least $12$. \\

Now we count the norm $12$ lattice vectors in this lattice ${\bf L}_{2n}$, it is clear there are ${n \choose 2} \cdot 9 \cdot 2$ type 1) norm $12$ lattice vectors in ${\bf L}_{2n}$. Each norm $12$ type 2) lattice vector in ${\bf L}_{2n}$ consists of six columns of the form $\pm(1,-1,0)^{\tau}$ or their permutations, the positions of these six non-zero columns correspond to the support set of a Hamming weight $6$ codeword in the ternary code ${\bf C}$. From condition 2) and 3) and the property that the sum of all coordinates of each codeword in ${\bf C}$ is zero, the sum of each row of any attice vector in ${\bf L}_{2n}$ can be  divided by $3$. Thus each Hamming weight $6$ codeword consists of six $\pm1$'s or three $+1$'s and three $-1$'s. Without loss of generality we treat the case of six $+1$'s. For each type 2)  norm $12$ lattice vector in ${\bf L}_{2n}$ such that the difference of the first row and the second row is a codeword in ${\bf C}$ consisting of six $1$'s, we count the number of such norm $12$ lattice vectors.\\

2A) There are two $\pm1$'s in the first row we have ${6 \choose 2} \cdot 2 =30$ such vectors. For example, if the first row is $(1,-1,0,0,0,0)$, we want the difference of the first row and the second row is $(1,1,1,1,1,1)$, then the second row has to be $(0,1,-1,-1,-1,-1)$. The same analysis for the first row $(-1,1,0,0,0,0)$.\\

2B) There are three $\pm1$'s in the first row we have ${6 \choose 3} \cdot 2 =40$ such vectors. For example, if the first row is $(1,1,1,0,0,0)$, we want the difference of the first row and the second row is $(1,1,1,1,1,1)$, then the second row has to be $(0,0,0,-1,-1,-1)$. The same analysis for the first row $(-1,-1,-1,0,0,0)$.\\

2C) There are four $\pm1$'s in the first row we have ${6 \choose 4} \cdot {4 \choose 2} =90$ such vectors. For example, if the first row is $(1,1,-1-1,0,0)$, we want the difference of the first row and the second row is $(1,1,1,1,1,1)$, then the second row has to be $(0,0,1,1,-1,-1)$. The same analysis for the first row $(-1,-1,1,1,0,0)$.\\

2D) There are five $\pm1$'s in the first row we have ${6 \choose 5} \cdot {5 \choose 1} \cdot 2 =60$ such vectors. For example, if the first row is $(1,1,1,1,-1,0)$, we want the difference of the first row and the second row is $(1,1,1,1,1,1)$, then the second row has to be  $(0,0,0,0,1,-1)$. The same analysis for the first row $(-1,-1,-1,-1,1,0)$.\\

2E) There are six $\pm1$'s in the first row we have ${6 \choose 3} \cdot+2=22$ such vectors. For example, if the first row is $(1,1,1,1,-1,-1,-1)$, we want the difference of the first row and the second row is $(1,1,1,1,1,1)$, then the second row has to be $(0,0,0,1,1,1)$. If the first row is $(1,1,1,1,1,1)$, the second row has to be $(0,0,0,0,0,0)$. The same analysis for the first row $(-1,-,1-1,1,1,1)$.\\

2F) There is no non-zero entry in the first row we have $2$ such vectors. The second  row has to be  $-(1,1,1,1,1,1)$, the third row has to be $\pm(1,1,1,1,1,1)$. We have one such norm $12$ type 2) lattice vectors.\\

Totally there are $30+40+90+60+22+1=243$ such type 2) norm $12$ lattice vectors for each Hamming weight $6$ codeword in ${\bf C}$. Thus the kissing number of the lattice ${\bf L}_{2n}$ is ${n \choose 2} \cdot 9 \cdot 2 +243 A_6$.\\

When $n=6$ and the  $[6,6]_3$ code is used, we construct a $12$ dimensional lattice ${\bf L}_{12}$ with the center density $\frac{3^6}{3 \cdot 3^5 \cdot 3^3}=\frac{1}{27}$ and the kissing number $270+486=756$. This is the Coxeter-Todd lattice ${\bf K}_{12}$ agian. When $n=9$ and the $[9,3,6]_3$ code is used, the center density of the $18$ dimensional lattice ${\bf L}_{18}$ is $\frac{3^9}{3^{4.5} \cdot 3 \cdot 3^6}=\frac{1}{9\sqrt{3}}$ and the kissing number of ${\bf L}_{18}$ is $648+243 \cdot 24=6480$. This is the Kappa lattice ${\bf K}_{18}$ (\cite{Nebe}). When $n=10$ and the $[10,4,6]_3$ code is used,  the center density of the $20$ dimensional lattice ${\bf L}_{20}$ is $\frac{3^{10}}{3 \cdot 3^6 \cdot 3^5}=\frac{1}{9}$ and the kissing number of ${\bf L}_{20}$ is $810+243 \cdot 60=15390$. This is the Kappa lattice ${\bf K}_{20}$ (\cite{Nebe}).\\

When $n=11$ and the $[11,5,6]_3$ code is ued, we construct a $22$ dimensional lattice with the center density $\frac{3^{11}}{3 \cdot 3^6 \cdot 3^{5.5}}=\frac{1}{3\sqrt{3}}$ and the kissing number $990+243 \cdot 132=33066$.\\

{\bf Remark.} Actually the condition that the sum of all coordinates of the ternary code ${\bf C}$ is zero is not necessary. The resulted lattice has the same kissing number and the same volume. We guess these lattices are isometry to the lattices constructed in Theorem 11.1.\\

\section{New $40$ dimensional extremal even unimodular lattices}

{\bf Proposition 11.1.} {\em Suppose there is a self-dual ternary $[n, \frac{n}{2}, 6]_3$ code ${\bf C}$. Let ${\bf W}_n$ be the lattice from theorem 2.2 by using the ternary code ${\bf C}$ and the lattice ${\bf D}_n$. The $\frac{1}{3}{\bf W}_n$ is an integral even lattice with the center density $\frac{1}{4}$ and minimum norm $4$.}\\

The proof is direct.\\

It is well-known that there are six $[20, 10, 6]_3$ self-dual ternary codes with the weight distribution $A_6=120, A_9=4360, A_{12}=26280, A_{15}=25728$ and $A_{18}=2560$ (\cite{PSW}). All these codes can be generated by the Hamming weight $6$ codewords (\cite{PSW}, check the generator matrices in the paper). We have an integral even lattice ${\bf W(C)}_{40}$ with the kissing number $k({\bf W}_{40})=20 \cdot 19 \cdot 2 \cdot 3+120 \cdot 183=24240$ from Proposition 8.1 and each ternary self-dual $[20, 10, 6]_3$ code ${\bf C}$ in \cite{PSW}. By using the following four translates of the lattice ${\bf W(C)}_{40}$, we get an unimodular even lattice ${\bf T(C)}_{40}$ of minimum norm $4$ naturally. Let ${\bf x}=(\frac{3}{2},\frac{3}{2},\frac{3}{2},\frac{3}{2},\frac{3}{2},\frac{3}{2},\frac{3}{2},\frac{3}{2},\frac{3}{2},\frac{3}{2}$,
$-\frac{3}{2},-\frac{3}{2},-\frac{3}{2},-\frac{3}{2},$ $-\frac{3}{2},-\frac{3}{2},-\frac{3}{2},-\frac{3}{2},-\frac{3}{2},-\frac{3}{2})$, and ${\bf X}_1=({\bf x}, -{\bf x},{\bf 0})$, ${\bf X}_2=(-{\bf x}$, ${\bf 0},{\bf x}), {\bf X}_3=({\bf 0}, -{\bf x},{\bf x})$. The four translates, ${\bf W(C)}_{40}, {\bf X}_1+{\bf W(C)}_{40}, {\bf X}_2+{\bf W(C)}_{40}$ and ${\bf X}_3+{\bf W(C)}_{40}$ is a lattice ${\bf T(C)}_{40}'$ with the volume $3^{40}$. The minimum norm of each vector in the translates ${\bf X}_i +{\bf W(C)}_{40}$ is at least $\frac{3}{2} \cdot 18 +\frac{9}{2} \cdot 2=36$. Thus the lattice ${\bf T(C)}_{40}=\frac{1}{3}{\bf T(C)}_{40}'$ is an unimodular even lattice with the minimum norm $4$. The kissing number of the lattice is $24240+3 \cdot2  \cdot 2560=39600$.\\

{\bf Proposition 11.2.} {\em The lattice ${\bf T(C)}_{40}$ is generated by the lattice vectors of minimum norm $36$. The automorphism group of the corresponding ternary self-dual $[20, 10, 6]_3$ code ${\bf C}$ is a subgroup of the automorphism group of the lattice ${\bf T(C)}_{40}$.}\\

{\bf Proof.} It is obvious the lattice ${\bf W(C)}_{40}$ is generated by the minimum norm vectors, since the ternary self-dual $[20, 10, 6]_3$ code is generated by the Hamming weight $6$ codewords. It is obvious that the $3$ vectors ${\bf X}_1,{\bf X}_2$ and ${\bf X}_3$ can be represented as the integral coefficient linear combinations of the minimum norm lattice vectors arising from the weight $18$ codewords and the minimum norm lattices vectors of the form $3 {\bf y}$. Then the first conclusion follows directly. The second conclusion is obvious.\\

{\bf Theorem 11.3.} {\em If two extremal even unimodular lattices of dimension $40$  ${\bf T(C)}_{40}$ and ${\bf T(C')}_{40}$ are isomorphic, then the two ternary self-dual $[20, 10, 6]_3$ codes ${\bf C}$ and ${\bf C'}$ are isomorphic.}\\

{\bf Proof.} Let ${\bf E}$ be the integral lattice ${\bf E}=\{({\bf x}_1,{\bf x}_2,{\bf x}_3):{\bf x}_i \in {\bf D}_{20}, {\bf x}_1+{\bf x}_2+{\bf x}_3=0\}$. If ${\bf G}$ is an orthogonal transformation satisfying ${\bf G}({\bf T(C)}_{40})={\bf T(C')}_{40}$, we prove that ${\bf G}(3{\bf E})=3{\bf E}$. Let $min({\bf T(C)}_{40}')$ be the set of norm $36$ lattice vectors of the lattice ${\bf T(C)}_{40}'$. There are $3 \frac{20 \cdot 19} {2} \cdot 4=2280$ type I such lattice vectors of the form $3 {\bf y}$ where ${\bf y}$ is an Euclid norm $4$ lattice vector in the lattice ${\bf E}$, $120 \cdot 183=21960$ type II such lattice vectors arising from the Hamming weight $6$ codewords and $ 2560 \cdot 2 \cdot 3=15360$ type III such lattice vectors in the three translates arising from the Hamming weight $18$ codewords. For each type III minimum norm lattice vector ${\bf U}_{III}$ in the lattice ${\bf T(C)}_{40}'$ corresponding to a weight $18$ codeword ${\bf c}$, we can check that there exist at most  $1+2+2N+t$ minimum norm lattice vectors ${\bf U}$ such that the sum ${\bf U}_{III}+{\bf U}$ is another minimum norm lattice vector, where $N$ is the number of weight $6$ codewords  ${\bf c}'$ satisfying the property the support $Supp({\bf c}')$  is not in the support $Supp({\bf c})$ of the codeword ${\bf c}$ and $wt({\bf c}-{\bf c}')$ is a weight $18$ codeword. Here $t$ is the number of the number of Hamming weight $6$ codewords such that $wt({\bf c}+{\bf c}')$ is a weight $18$ codeword and $Supp({\bf c}') \subset Supp({\bf c})$. Among these $1+2+2N+t$ minimum norm lattice vectors, one is a type I minimum norm lattice vector, $N+t$'s are type II minimum norm lattice vectors and the other $2+N$ are type III minimum norm lattice vectors. From the paper $\cite{PSW}$ we have $N \leq 36$ and $t \leq 24$. On the other hand for any type I or II minimum norm lattice vector, there are  more than $120$ minimum norm lattice vectors satisfying this property. \\

Thus the orthogonal transformation ${\bf G}$ has to send a type III minimum norm lattice vector of the lattice ${\bf T(C)}_{40}'$ to a type III minimum norm lattice vector in the lattice ${\bf T(C')}_{40}'$. Moreover ${\bf G}$ sends the corresponding type I minimum norm lattice vector of the lattice ${\bf T(C)}_{40}'$ to a type I minimum norm lattice vector of the lattice ${\bf T(C')}_{40}'$, since for any one of these $N$ type II minimum lattice vector ${\bf U}_{II}$, there are more type III minimum norm lattice vector ${\bf U}_{III}$ such that the sum ${\bf U}_{II}+{\bf U}_{III}$ is another minimum norm lattice vector. This implies that ${\bf G}$ always sends the sub-lattice $3{\bf E}$ of the lattice ${\bf T(C)}_{40}'$ to the sub-lattice $3{\bf E}$ of the lattice ${\bf T(C')}_{40}'$. Then ${\bf G}$ induces an isomorphism of the code ${\bf C}$ to the code ${\bf C'}$ naturally.\\

Since the Mckay extremal even unimodular lattice (\cite{CS1}, page 221) and the Ozeki extremal even unimodular lattices
(\cite{Ozeki}) of dimension $40$ are not generated by the minimum norm lattice vectors, the extremal even unimodular lattices in Theorem 3.3
are not the same as Mckey and Ozeki lattices. The only two known extremal even unimodular lattices of dimension $40$ which are generated
by the minimum norm lattice vectors are the lattice of Calderbank and Sloane (\cite{Calderbank}) and the lattice of G. Nebe (\cite{Nebe}). From Theorem 11.3 at least four new $40$ dimensional
extremal even unimodular lattices are constructed.\\

{\bf Remarkn 11.1.} If we set the vector ${\bf x}=(-\frac{3}{2}, \frac{3}{2},...,\frac{3}{2}) \in {\bf R}^{20}$, it seems new $40$ dimensional extremal even unimodular lattices can be constructed. \\

\section{New $32$ dimensional extremal even unimodular lattices}

In this section we give four constructed even unimodular lattices of dimensions $32$, however we do not know how to decide if they are distinct or they are different with the 15 Koch-Venkov extremal even unimodular lattices of dimension $32$(\cite{KochVen,KochNebe}. We just use Magma to help us to decide the isometric problem.\\

{\bf 1st}.\\

In our ternary construction we use the following $32$ dimensional lattice ${\bf L} \subset {\bf Z}^{16} \oplus {\bf Z}^{16} \oplus {\bf Z}^{16}$ defined by $x_{1j}+x_{2j}+x_{3j}=0$ for $j=1,...,16$. This is a lattice with minimum norm $2$. A ternary self-dual $[16, 8, 3]_3$ code ${\bf C}$ is used in the teranry construction. This ${\bf C}$ is the direct sum of the linear self-dual $[4,2,3]_3$ code with the following generator matrix. Four copies pf such ternary linear code are supported at the four position sets $\{1,2,3,4\}$, $\{5,6,7,8\}$, $\{9,10,11,12\}$ and $\{13,14,15,16\}$.\\

$$
\left(
\begin{array}{cccc}
1&1&1&0\\
0&1&-1&1\\
\end{array}
\right)
$$

We get a $32$ dimensional lattice ${\bf T'}_{32}$ with minimum norm $18$ and volume $3^{32}$. The lattice is the union of $3^8$ translates of the lattice $3{\bf L}$ with leading vectors from the $3^8$ codewords in the ternary code ${\bf C}$.\\

We need to use $2^{6+10}$ translates of this lattice ${\bf T'}_{32}$ according to the following two binary codes ${\bf C}_1$ and ${\bf C}_2$. Here ${\bf C}_1$ is the binary $[16, 6, 4]_2$ code with the following generator matrix.\\

$$
\left(
\begin{array}{cccccccccccccccccc}
1&1&1&1&1&1&1&1&0&0&0&0 &0&0&0&0\\
0&0&0&0 &1&1&1&1 &1&1&1&1& 0&0&0&0\\
1&0&0&0& 1&0&0&0& 1&0&0&0& 1&0&0&0\\
0&1&0&0& 0&1&0&0& 0&1&0&0& 0&1&0&0\\
0&0&1&0& 0&0&1&0& 0&0&1&0& 0&0&1&0\\
0&0&0&1& 0&0&0&1& 0&0&0&1& 0&0&0&1\\
\end{array}
\right)
$$
 
The code ${\bf C}_2$ is a binary $[16,10,4]_2$ code with the following generator matrix.\\

$$
\left(
\begin{array}{cccccccccccccccccc}
1&1&1&1&1&1&1&1&0&0&0&0 &0&0&0&0\\
0&0&0&0 &1&1&1&1 &1&1&1&1& 0&0&0&0\\
1&0&0&0& 1&0&0&0& 1&0&0&0& 1&0&0&0\\
0&1&0&0& 0&1&0&0& 0&1&0&0& 0&1&0&0\\
0&0&1&0& 0&0&1&0& 0&0&1&0& 0&0&1&0\\
0&0&0&1& 0&0&0&1& 0&0&0&1& 0&0&0&1\\
1&1&0&0& 1&1&0&0& 0&0&0&0& 0&0&0&0\\
1&1&0&0& 0&0&0&0& 1&1&0&0& 0&0&0&0\\
0&1&1&0& 0&1&1&0& 0&0&0&0& 0&0&0&0\\
0&1&1&0& 0&0&0&0& 0&1&1&0& 0&0&0&0\\
\end{array}
\right)
$$

We note that the binary linear code ${\bf C}_2$ is generated by the vector $(1,0,0,0,1,0,0,0,1,0,0,0,1,0,0,0) \in {\bf F}_2^{16}$ and the binary linear $[16, 9, 4]$ code ${\bf C'}_2$ with the following parity check matrix.\\

$$
\left(
\begin{array}{cccccccccccccccccc}
1&0&0&0&1&0&0&0&1&0&0&0 &1&0&0&0\\
0&1&0&0 &0&1&0&0 &0&1&0&0& 0&1&0&0\\
0&0&1&0& 0&0&1&0& 0&0&1&0& 0&0&1&0\\
0&0&0&1& 0&0&0&1& 0&0&0&1& 0&0&0&1\\
1&1&1&1& 0&0&0&0& 0&0&0&0& 0&0&0&0\\
0&0&0&0& 1&1&1&1& 0&0&0&0& 0&0&0&0\\
0&0&0&0& 0&0&0&0& 1&1&1&1& 0&0&0&0\\
\end{array}
\right)
$$

That is, the sums of coordinates of any codeword in ${\bf C'}_2$ at the position sets $\{1,2,3,4\}, \{5,6,7,8\}$, $\{9,10,11,12\}$ and $\{13,14,15,16\}$ are zero, the sums of coordinates of any codeword in ${\bf C'}_2$ at the position sets $\{1,5,9,13\}$. $\{2,6,10.14\}$, $\{3,7,11,15\}$ and $\{4, 8, 12,16\}$ are zero.\\

The union of these $2^{6+10}$ translates is a new $32$ dimensional lattice ${\bf T}_{32}$ with the volume $\frac{3^{32}}{2^{18}}$.   For any codeword $c \in {\bf C}_1$, whose support positions are $i_1,...,i_t$, set ${\bf x}_c$ the vector in ${\bf R}^{16}$ with $\frac{3}{2}$ at these positions $i_1,..,i_t$. Let ${\bf X}_{c, 1}=({\bf x}_c, -{\bf x}_c,{\bf 0})$.  For any codeword $c \in {\bf C}_2$, whose support positions are $i_1,...,i_t$, set ${\bf x}_c$ the vector in ${\bf R}^{16}$ with $\frac{3}{2}$ at these positions $i_1,..,i_t$. Let ${\bf X}_{c, 2}=({\bf 0}, {\bf x}_c, -{\bf x}_c)$.  Totally we get $2^{10+6}=2^{16}$ leading vectors in the real space spanned by the lattice ${\bf T'}_{32}$. The lattice ${\bf T}_{32}$ is the union of $2^{16}$ translates of ${\bf T'}_{32}$ with these $2^{16}$ leading vectors. Its volume is $\frac{3^{32}}{2^{16}}$\\

{\bf Theorem 13.1.} {\em The norms of the differences in the above $2^{16}$ translates are at least $18$. The center density of this lattice ${\bf T}_{32}$ is $1$.}\\

{\bf Proof.} For any ${\bf X}_{c,1}-{\bf X}_{c',2}-{\bf e}$, where  $e$ is a lattice vector in ${\bf T'}_{32}$, if ${\bf e}$ is in $3{\bf L}$, the norm of this difference is at least $18$.  If ${\bf e}$ is not in $3{\bf L}$,  it is from a cordword  of weights  $6, 9, 12$ in the ternary code ${\bf C}$.\\

It is clear that every non-zero codeword in the ternary $[4,2,3]_3$ code described as above is with Hamming weight $3$. Thus the minimum norms of the differences ${\bf X}_{c,1}-{\bf X}_{c',2}-{\bf e}$ at each of the four position sets $\{1,2,3,4\}$, $\{5,6,7,8\}$, $\{9,10,11,12\}$ and $\{13,14,15,16\}$ are at least $\frac{3}{2} \cdot 3=\frac{9}{2}$. This minimum norm possibility $\frac{9}{2}$ happens only if the union of the supports of $c$ and $c'$ at this position set has three element. If this minimum norm possibility does not happen, that is, the union of the supports of the codewords $c$ and $c'$ in this four-element position set has one or two or four elements, the norm of the difference is at least $9$.\\

We observe that one codeword in the linear binary code ${\bf C}_1$ is either not zero at each of the above four position sets, or is non-zero at all four elements of two of the above four position sets. Thus if $c$ is non-zero, the minimum norm of the difference is at least $9+9=\frac{9}{2} \cdot 4=18$. If $c$ is zero and $c'$ is not zero, we observe that if $c'$ is in the linear binary code ${\bf C'}_2$, the minimum norm of the difference is at least $9+9$ (that is, the support of $c'$ in each of the above four position sets has two or four elements). If $c'$ is not in the linear binary code ${\bf C'}_2$, the vector $(1,0,0,0,1,0,0,0,1,0,0,0,1,0,0,0)$ is involved the minimum norm of the difference is at least $\frac{9}{2} \cdot 4=18$.\\

{\bf Corollary 13.2} {\em The lattice $\frac{\sqrt{2}}{3} {\bf T}_{32}$ is an extremal even uni-modular lattice of dimension $32$.}\\

The following is a base of this lattice.\\

v1 := V![3/2,3/2,3/2,3/2, 3/2,3/2,3/2,3/2, 3/2,3/2,3/2,3/2, 3/2,3/2,3/2,3/2, -3/2,-3/2,-3/2,-3/2, -3/2,-3/2,-3/2,-3/2, -3/2,-3/2,-3/2,-3/2, 
-3/2,-3/2,-3/2,-3/2, 0,0,0,0, 0,0,0,0, 0,0,0,0, 0,0,0,0];\\

v2 := V![3/2,0,0,0, 3/2,0,0,0, 3/2,0,0,0, 3/2,0,0,0, -3/2,0,0,0, -3/2,0,0,0, -3/2,0,0,0, -3/2,0,0,0, 0,0,0,0, 
0,0,0,0, 0,0,0,0, 0,0,0,0];\\

v3 := V![0, 3/2,0,0, 0,3/2,0,0, 0,3/2,0,0, 0,3/2,0,0, 0,-3/2,0,0, 0,-3/2,0,0, 0,-3/2,0,0, 0,-3/2,0,0, 0,0,0,0,
0,0,0,0, 0,0,0,0, 0,0,0,0];\\

v4 := V![0,0,3/2,0, 0,0,3/2,0, 0,0,3/2,0, 0,0,3/2,0, 0,0,-3/2,0, 0,0,-3/2,0, 0,0,-3/2,0, 0,0,-3/2,0, 0,0,0,0,
0,0,0,0, 0,0,0,0, 0,0,0,0];\\

v5 := V![3/2,3/2,3/2,3/2, 3/2,3/2,3/2,3/2, 0,0,0,0, 0,0,0,0, -3/2,-3/2,-3/2,-3/2, -3/2,-3/2,-3/2,-3/2, 0,0,0,0, 0,0,0,0, 
0,0,0,0, 0,0,0,0, 0,0,0,0, 0,0,0,0];\\

v6 := V![0,0,0,0, 3/2,3/2,3/2,3/2, 3/2,3/2,3/2,3/2, 0,0,0,0, 0,0,0,0, -3/2,-3/2,-3/2,-3/2, -3/2,-3/2,-3/2,-3/2, 0,0,0,0, 
0,0,0,0, 0,0,0,0, 0,0,0,0, 0,0,0,0];\\

v7 := V![3/2,3/2,0,0, 3/2,3/2,0,0,  0,0,0,0, 0,0,0,0, -3/2,-3/2,0,0,  -3/2,-3/2,0,0, 0,0,0,0, 0,0,0,0, 0,0,0,0, 0,0,0,0,
0,0,0,0, 0,0,0,0];\\

v8 := V![3/2,3/2,0,0, 0,0,0,0, 3/2,3/2,0,0,  0,0,0,0, -3/2,-3/2,0,0, 0,0,0,0,  -3/2,-3/2,0,0, 0,0,0,0, 0,0,0,0, 0,0,0,0, 
0,0,0,0, 0,0,0,0];\\

v9 := V![0,3/2,3/2,0, 0,3/2,3/2,0, 0,0,0,0, 0,0,0,0, 0,-3/2,-3/2,0, 0,-3/2,-3/2,0, 0,0,0,0, 0,0,0,0, 0,0,0,0,
0,0,0,0, 0,0,0,0, 0,0,0,0];\\

v10 := V![0,3/2,3/2,0, 0,0,0,0, 0,3/2,3/2,0, 0,0,0,0, 0,-3/2,-3/2,0, 0,0,0,0, 0,-3/2,-3/2,0, 0,0,0,0, 0,0,0,0, 0,0,0,0,
0,0,0,0, 0,0,0,0];\\

v11 := V![0,0,0,0, 0,0,0,0, 0,0,0,0, 0,0,0,0, 3/2,3/2,3/2,3/2, 3/2,3/2,3/2,3/2, 3/2,3/2,3/2,3/2, 3/2,3/2,3/2,3/2, -3/2,-3/2,-3/2,-3/2, -3/2,-3/2,-3/2,-3/2, -3/2,-3/2,-3/2,-3/2, 
-3/2,-3/2,-3/2,-3/2];\\

v12 := V![0,0,0,0, 0,0,0,0, 0,0,0,0, 0,0,0,0, 3/2,0,0,0, 3/2,0,0,0, 3/2,0,0,0, 3/2,0,0,0,
 -3/2,0,0,0, -3/2,0,0,0, -3/2,0,0,0, -3/2,0,0,0];\\

v13 := V![0,0,0,0, 0,0,0,0, 0,0,0,0, 0,0,0,0, 3/2,3/2,3/2,3/2, 3/2,3/2,3/2,3/2, 0,0,0,0,0,0,0,0, 
-3/2,-3/2,-3/2,-3/2, -3/2,-3/2,-3/2,-3/2, 0,0,0,0, 0,0,0,0];\\

v14 := V![ 0,0,0,0, 0,0,0,0, 0,0,0,0, 0,0,0,0,  0,0,0,0, 3/2,3/2,3/2,3/2, 3/2,3/2,3/2,3/2, 0,0,0,0,
0,0,0,0, -3/2,-3/2,-3/2,-3/2, -3/2,-3/2,-3/2,-3/2, 0,0,0,0];\\

v15 :=V![-1/2,-1/2,-1/2,0, -1/2,-1/2,-1/2,0, -1/2,-1/2,-1/2,0,-1/2,-1/2,-1/2,0, 1,1,1,0, 1,1,1,0, 1,1,1,0, 1,1,1,0, 
-1/2,-1/2,-1/2,0, -1/2,-1/2,-1/2,0, -1/2,-1/2,-1/2,0,-1/2,-1/2,-1/2,0];\\

v16 :=V![-1/2,1/2,0,-1/2, -1/2,1/2,0,-1/2, -1/2,1/2,0,-1/2, -1/2,1/2,0,-1/2, 1,-1,0,1, 1,-1,0,1, 1,-1,0,1, 1,-1,0,1, 
-1/2,1/2,0,-1/2, -1/2,1/2,0,-1/2, -1/2,1/2,0,-1/2, -1/2,1/2,0,-1/2];\\

v17 :=V![1,1,1,0, 0,0,0,0, 0,0,0,0, 0,0,0,0, -2,-2,-2,0, 0,0,0,0, 0,0,0,0, 0,0,0,0, 1,1,1,0, 0,0,0,0, 0,0,0,0, 0,0,0,0];\\

v18 :=V![1,-1,0,1, 0,0,0,0, 0,0,0,0, 0,0,0,0, -2,2,0,-2,0,0,0,0, 0,0,0,0, 0,0,0,0, 1,-1,0,1, 0,0,0,0, 0,0,0,0, 0,0,0,0];\\

v19 :=V![0,0,0,0, 1,1,1,0, 0,0,0,0, 0,0,0,0,  0,0,0,0, -2,-2,-2,0, 0,0,0,0, 0,0,0,0,  0,0,0,0, 1,1,1,0, 0,0,0,0, 0,0,0,0];\\

v20 :=V![0,0,0,0, 1,-1,0,1, 0,0,0,0, 0,0,0,0, 0,0,0,0, -2,2,0,-2,0,0,0,0, 0,0,0,0, 0,0,0,0, 1,-1,0,1, 0,0,0,0, 0,0,0,0];\\

v21 :=V![0,0,0,0, 0,0,0,0, 1,1,1,0, 0,0,0,0, 0,0,0,0, 0,0,0,0, -2,-2,-2,0, 0,0,0,0, 0,0,0,0, 0,0,0,0, 1,1,1,0, 0,0,0,0];\\

v22 :=V![0,0,0,0, 0,0,0,0, 1,-1,0,1, 0,0,0,0, 0,0,0,0, 0,0,0,0, -2,2,0,-2, 0,0,0,0, 0,0,0,0, 0,0,0,0, 1,-1,0,1, 0,0,0,0];\\

v23 :=V![3,0,0,0, 0,0,0,0, 0,0,0,0, 0,0,0,0, -3,0,0,0, 0,0,0,0, 0,0,0,0, 0,0,0,0, 0,0,0,0, 0,0,0,0, 0,0,0,0, 0,0,0,0];\\

v24 :=V![0,3,0,0, 0,0,0,0, 0,0,0,0, 0,0,0,0, 0,-3,0,0, 0,0,0,0, 0,0,0,0, 0,0,0,0, 0,0,0,0, 0,0,0,0, 0,0,0,0, 0,0,0,0];\\

v25 :=V![0,0,3,0, 0,0,0,0, 0,0,0,0, 0,0,0,0, 0,0,-3,0, 0,0,0,0, 0,0,0,0, 0,0,0,0, 0,0,0,0, 0,0,0,0, 0,0,0,0, 0,0, 0,0];\\

v26 :=V![0,0,0,3, 0,0,0,0, 0,0,0,0, 0,0,0,0, 0,0,0,-3, 0,0,0,0, 0,0,0,0, 0,0,0,0, 0,0,0,0, 0,0,0,0, 0,0,0,0, 0,0,0,0];\\

v27 :=V![0,0,0,0, 3,0,0,0, 0,0,0,0, 0,0,0,0, 0,0,0,0, -3,0,0,0, 0,0,0,0, 0,0,0,0, 0,0,0,0, 0,0,0,0, 0,0,0,0, 0,0,0,0];\\

v28 :=V![0,0,0,0, 0,0,0,0, 3,0,0,0, 0,0,0,0, 0,0,0,0, 0,0,0,0, -3,0,0,0, 0,0,0,0, 0,0,0,0, 0,0,0,0, 0,0,0,0, 0,0,0,0];\\

v29 :=V![0,0,0,0, 0,0,0,0, 0,0,0,0, 0,0,0,0, 3,0,0,0, 0,0,0,0, 0,0,0,0, 0,0,0,0, -3,0,0,0, 0,0,0,0, 0,0,0,0, 0,0,0,0];\\

v30 :=V![0,0,0,0, 0,0,0,0, 0,0,0,0, 0,0,0,0, 0,3,0,0, 0,0,0,0, 0,0,0,0, 0,0,0,0, 0,-3,0,0, 0,0,0,0, 0,0,0,0, 0,0,0,0];\\

v31 :=V![0,0,0,0, 0,0,0,0, 0,0,0,0, 0,0,0,0, 0,0,0,0, 3,0,0,0, 0,0,0,0, 0,0,0,0, 0,0,0,0, -3,0,0,0, 0,0,0,0, 0,0,0,0];\\

v32 :=V![0,0,0,0, 0,0,0,0, 0,0,0,0, 0,0,0,0,  0,0,0,0, 0,0,0,0, 3,0,0,0, 0,0,0,0,  0,0,0,0, 0,0,0,0, -3,0,0,0, 0,0,0,0];\\

From Magma we can compute  the Gram matrix of the lattice $\frac{\sqrt{2}}{9} {\bf T}_{32}$ and its kissing number is indeed $146880$. It is easy to see that $3^2$ is a factor of the order of the automorphism group of this lattice.\\

{\bf Another form}\\

In our ternary construction we use the following $32$ dimensional lattice ${\bf L} \subset {\bf Z}^{16} \oplus {\bf Z}^{16} \oplus {\bf Z}^{16}$ defined by $x_{1j}+x_{2j}+x_{3j}=0$ for $j=1,...,16$. This is a lattice with minimum norm $2$. A ternary self-dual $[16, 8, 3]_3$ code ${\bf C}$ is used in the teranry construction. This ${\bf C}$ is the direct sum of tfour copies of the linear self-dual $[4,2,3]_3$ code. \\

We get a $32$ dimensional lattice ${\bf T'}_{32}^2$ with minimum norm $18$ and volume $3^{32}$. The lattice is the union of $3^8$ translates of the lattice $3{\bf L}$ with leading vectors from the $3^8$ codewords in the ternary code ${\bf C}$.\\

We need to use $2^{16}$ translates of this lattice ${\bf T'}_{32}^2$ according to the following linear binary code ${\bf C}$. Here ${\bf C}_3$ is the binary $[16, 8, 4]_2$ code with the following generator matrix.\\

$$
\left(
\begin{array}{cccccccccccccccccc}
1&1&0&0&1&1&0&0&0&0&0&0 &0&0&0&0\\
0&0&0&0 &1&1&0&0 &1&1&0&0& 0&0&0&0\\
0&0&0&0&0&0&0 &0&1&1&0&0 &1&1&0&0\\
0&0&1&1 &0&0&1&1 &0&0&0&0& 0&0&0&0\\
0&0&0&0 &0&0&1&1 &0&0&1&1& 0&0&0&0\\
0&0&0&0 &0&0&0&0 &0&0&1&1& 0&0&1&1\\
1&0&0&0& 1&0&0&0& 1&0&0&0& 1&0&0&0\\
0&0&1&0& 0&0&1&0& 0&0&1&0& 0&0&1&0\\
\end{array}
\right)
$$

The union of these $2^{8+8}$ translates is a new $32$ dimensional lattice ${\bf T}_{32}^2$ with the volume $\frac{3^{32}}{2^{16}}$.   For any codeword $c \in {\bf C}_3$, whose support positions are $i_1,...,i_t$, set ${\bf x}_c$ the vector in ${\bf R}^{16}$ with $\frac{3}{2}$ at these positions $i_1,..,i_t$. Let ${\bf X}_{c, 1}=({\bf x}_c, -{\bf x}_c,{\bf 0})$.  For any codeword $c \in {\bf C}_3$, whose support positions are $i_1,...,i_t$, set ${\bf x}_c$ the vector in ${\bf R}^{16}$ with $\frac{3}{2}$ at these positions $i_1,..,i_t$. Let ${\bf X}_{c, 2}=({\bf 0}, {\bf x}_c, -{\bf x}_c)$.  Totally we get $2^{10+6}=2^{16}$ leading vectors in the real space spanned by the lattice ${\bf T'}_{32}$. The lattice ${\bf T}_{32}^2$ is the union of $2^{16}$ translates of ${\bf T'}_{32}$ with these $2^{16}$ leading vectors. Its volume is $\frac{3^{32}}{2^{16}}$\\

{\bf Theorem 13.3.} {\em The norms of the differences in the above $2^{16}$ translates are at least $18$. The center density of this lattice ${\bf T}_{32}^2$ is $1$.}\\

{\bf Proof.} For any ${\bf X}_{c,1}-{\bf X}_{c',2}-{\bf e}$, where  $e$ is a lattice vector in ${\bf T'}_{32}$, if ${\bf e}$ is in $3{\bf L}$, the norm of this difference is at least $18$.  If ${\bf e}$ is not in $3{\bf L}$,  it is from a cordword  of weights  $6, 9, 12$ in the ternary code ${\bf C}$.\\

It is clear that every non-zero codeword in the ternary $[4,2,3]_3$ code described as above is with Hamming weight $3$. Thus the minimum norms of the differences ${\bf X}_{c,1}-{\bf X}_{c',2}-{\bf e}$ at each of the four position sets $\{1,2,3,4\}$, $\{5,6,7,8\}$, $\{9,10,11,12\}$ and $\{13,14,15,16\}$ are at least $\frac{3}{2} \cdot 3=\frac{9}{2}$. This minimum norm possibility $\frac{9}{2}$ happens only if the union of the supports of $c$ and $c'$ at this position set has three element. If this minimum norm possibility does not happen, that is, the union of the supports of the codewords $c$ and $c'$ in this four-element position set has one or two or four elements, the norm of the difference is at least $9$.\\

We observe that one codeword in the linear binary code ${\bf C}_3$ is either not zero at each of the above four position sets, or only have two non-zero positions at  two of the above four-element-position sets. Thus if $c$ is non-zero at the all four-element-position sets, the minimum norm of the difference is at least $\frac{9}{2} \cdot 4=18$. If it only have two non-zero positions at two of the four-element-position sets, the minimum norm of the difference is at least $9+9=18$. The conlusion is proved.\\

{\bf Corollary 13.4.} {\em $\frac{\sqrt{2}}{3}{\bf T}_{32}^2$ is an unimodular lattice of dimension $32$. Its minimum norm is $4$.}\\

The following is a base of this lattice. It is easy to see that $3^3$ is a factor of the order of the automorphism of this lattice.\\

v1 := V![3/2,3/2,0,0, 3/2,3/2,0,0, 0,0,0,0, 0,0,0,0, -3/2,-3/2,0,0, -3/2,-3/2,0,0, 0,0,0,0,
0,0,0,0, 0,0,0,0, 0,0,0,0, 0,0,0,0, 0,0,0,0];\\

v2 := V![0,0,3/2,3/2, 0,0,3/2,3/2, 0,0,0,0, 0,0,0,0,  0,0, -3/2,-3/2, 0,0,-3/2,-3/2, 0,0,0,0, 0,0,0,0,
0,0,0,0, 0,0,0,0, 0,0,0,0, 0,0,0,0];\\

v3 := V![0,0,0,0, 3/2,3/2,0,0, 3/2,3/2,0,0, 0,0,0,0, 0,0,0,0, -3/2,-3/2,0,0, -3/2,-3/2,0,0, 0,0,0,0, 0,0,0,0,
0,0,0,0, 0,0,0,0, 0,0,0,0];\\

v4 := V![0,0,0,0, 0,0,3/2,3/2, 0,0,3/2,3/2, 0,0,0,0, 0,0,0,0, 0,0,-3/2,-3/2, 0,0,-3/2,-3/2, 0,0,0,0, 0,0,0,0,
0,0,0,0, 0,0,0,0, 0,0,0,0];\\

v5 := V![0,0,0,0, 0,0,0,0, 3/2,3/2,0,0, 3/2,3/2,0,0, 0,0,0,0, 0,0,0,0, -3/2,-3/2,0,0, -3/2,-3/2,0,0,
0,0,0,0, 0,0,0,0, 0,0,0,0, 0,0,0,0];\\

v6 := V![0,0,0,0, 0,0,0,0, 0,0,3/2,3/2, 0,0,3/2,3/2, 0,0,0,0, 0,0,0,0, 0,0,-3/2,-3/2, 0,0,-3/2,-3/2, 
0,0,0,0, 0,0,0,0, 0,0,0,0, 0,0,0,0];\\

v7 := V![3/2,0,0,0, 3/2,0,0,0,  3/2,0,0,0, 3/2,0,0,0, -3/2,0,0,0,  -3/2,0,0,0, -3/2,0,0,0, -3/2,0,0,0, 0,0,0,0, 0,0,0,0,
0,0,0,0, 0,0,0,0];\\

v8 := V![0,0,3/2,0, 0,0,3/2,0, 0,0,3/2,0, 0,0,3/2,0, 0,0,-3/2,0, 0,0,-3/2,0, 0,0,-3/2,0, 0,0,-3/2,0, 
0,0,0,0, 0,0,0,0, 0,0,0,0, 0,0,0,0];\\

v9 := V![0,0,0,0, 0,0,0,0, 0,0,0,0, 0,0,0,0, 3/2,3/2,0,0, 3/2,3/2,0,0, 0,0,0,0, 0,0,0,0, -3/2,-3/2,0,0, -3/2,-3/2,0,0, 0,0,0,0, 0,0,0,0];\\

v10 := V![0,0,0,0, 0,0,0,0, 0,0,0,0, 0,0,0,0,  0,0,3/2,3/2, 0,0,3/2,3/2, 0,0,0,0, 0,0,0,0,  0,0,-3/2,-3/2,  0,0,-3/2,-3/2,  0,0,0,0, 0,0, 0,0];\\

v11 := V![0,0,0,0, 0,0,0,0, 0,0,0,0, 0,0,0,0, 0,0,0,0, 3/2,3/2,0,0, 3/2,3/2,0,0, 0,0,0,0, 0,0,0,0, -3/2,-3/2,0,0, -3/2,-3/2,0,0, 0,0,0,0];\\

v12 := V![0,0,0,0, 0,0,0,0, 0,0,0,0, 0,0,0,0,  0,0,0,0, 0,0,3/2,3/2, 0,0,3/2,3/2, 0,0,0,0,  0,0,0,0, 0,0,-3/2,-3/2, 0,0,-3/2,-3/2, 0,0,0,0];\\

v13 := V![0,0,0,0, 0,0,0,0, 0,0,0,0, 0,0,0,0, 0,0,0,0, 0,0,0,0, 3/2,3/2,0,0, 3/2,3/2,0,0, 0,0,0,0, 0,0,0,0, -3/2,-3/2,0,0, -3/2,-3/2,0,0];\\

v14 := V![0,0,0,0, 0,0,0,0, 0,0,0,0, 0,0,0,0, 0,0,0,0, 0,0,0,0, 0,0,3/2,3/2, 0,0,3/2,3/2, 0,0,0,0, 0,0,0,0, 0,0,-3/2,-3/2, 0,0,-3/2,-3/2];\\

v15 :=V![-1/2,-1/2,-1/2,0, -1/2,-1/2,-1/2,0, -1/2,-1/2,-1/2,0,-1/2,-1/2,-1/2,0, 1,1,1,0, 1,1,1,0, 1,1,1,0, 1,1,1,0, 
-1/2,-1/2,-1/2,0, -1/2,-1/2,-1/2,0, -1/2,-1/2,-1/2,0,-1/2,-1/2,-1/2,0];\\

v16 :=V![-1/2,0,1/2,1/2, -1/2,0,1/2,1/2, -1/2,0,1/2,1/2, -1/2,0,1/2,1/2, 1,0,-1,-1, 1,0,-1,-1, 1,0,-1,-1, 1,0,-1,-1, 
-1/2,0,1/2,1/2, -1/2,0,1/2,1/2, -1/2, 0, 1/2,   1/2, -1/2,0,1/2,1/2];\\

v17 :=V![1,1,1,0, 0,0,0,0, 0,0,0,0, 0,0,0,0, -2,-2,-2,0, 0,0,0,0, 0,0,0,0, 0,0,0,0, 1,1,1,0, 0,0,0,0, 0,0,0,0, 0,0,0,0];\\

v18 :=V![1,0,-1,-1, 0,0,0,0, 0,0,0,0, 0,0,0,0, -2,0,2,2,0,0,0,0, 0,0,0,0, 0,0,0,0, 1,0,-1,-1, 0,0,0,0, 0,0,0,0, 0,0,0,0];\\

v19 :=V![0,0,0,0, 1,1,1,0, 0,0,0,0, 0,0,0,0,  0,0,0,0, -2,-2,-2,0, 0,0,0,0, 0,0,0,0,  0,0,0,0, 1,1,1,0, 0,0,0,0, 0,0,0,0];\\

v20 :=V![0,0,0,0, 1,0,-1,-1, 0,0,0,0, 0,0,0,0, 0,0,0,0, -2,0,2,2,0,0,0,0, 0,0,0,0, 0,0,0,0, 1,0,-1,-1, 0,0,0,0, 0,0,0,0];\\

v21 :=V![0,0,0,0, 0,0,0,0, 1,1,1,0, 0,0,0,0, 0,0,0,0, 0,0,0,0, -2,-2,-2,0, 0,0,0,0, 0,0,0,0, 0,0,0,0, 1,1,1,0, 0,0,0,0];\\

v22 :=V![0,0,0,0, 0,0,0,0, 1,0,-1,-1, 0,0,0,0, 0,0,0,0, 0,0,0,0, -2,0,2,2, 0,0,0,0, 0,0,0,0, 0,0,0,0, 1,0,-1,-1, 0,0,0,0];\\

v23 :=V![3,0,0,0, 0,0,0,0, 0,0,0,0, 0,0,0,0, -3,0,0,0, 0,0,0,0, 0,0,0,0, 0,0,0,0, 0,0,0,0, 0,0,0,0, 0,0,0,0, 0,0,0,0];\\

v24 :=V![0,3,0,0, 0,0,0,0, 0,0,0,0, 0,0,0,0, 0,-3,0,0, 0,0,0,0, 0,0,0,0, 0,0,0,0, 0,0,0,0, 0,0,0,0, 0,0,0,0, 0,0,0,0];\\
v25 :=V![0,0,3,0, 0,0,0,0, 0,0,0,0, 0,0,0,0, 0,0,-3,0, 0,0,0,0, 0,0,0,0, 0,0,0,0, 0,0,0,0, 0,0,0,0, 0,0,0,0, 0,0, 0,0];\\

v26 :=V![0,0,0,3, 0,0,0,0, 0,0,0,0, 0,0,0,0, 0,0,0,-3, 0,0,0,0, 0,0,0,0, 0,0,0,0, 0,0,0,0, 0,0,0,0, 0,0,0,0, 0,0,0,0];\\

v27 :=V![0,0,0,0, 3,0,0,0, 0,0,0,0, 0,0,0,0,  0,0,0,0, -3,0,0,0, 0,0,0,0, 0,0,0,0, 0,0,0,0, 0,0,0,0, 0,0,0,0, 0,0,0,0];\\

v28 :=V![0,0,0,0, 0,0,3,0, 0,0,0,0, 0,0,0,0,  0,0,0,0, 0,0,-3,0, 0,0,0,0, 0,0,0,0, 0,0,0,0, 0,0,0,0, 0,0,0,0, 0,0,0,0];\\

v29 :=V![0,0,0,0, 0,0,0,0, 3,0,0,0, 0,0,0,0, 0,0,0,0, 0,0,0,0, -3,0,0,0, 0,0,0,0, 0,0,0,0, 0,0,0,0, 0,0,0,0, 0,0,0,0];\\

v30 :=V![0,0,0,0, 0,0,0,0, 0,0,3,0, 0,0,0,0, 0,0,0,0, 0,0,0,0, 0,0,-3,0, 0,0,0,0, 0,0,0,0, 0,0,0,0, 0,0,0,0, 0,0,0,0];\\

v31 :=V![0,0,0,0, 0,0,0,0, 0,0,0,0, 0,0,0,0, 3,0,0,0, 0,0,0,0, 0,0,0,0, 0,0,0,0, -3,0,0,0, 0,0,0,0, 0,0,0,0, 0,0,0,0];\\

v32 :=V![0,0,0,0, 0,0,0,0, 0,0,0,0, 0,0,0,0,  0,3,0,0, 0,0,0,0, 0,0,0,0, 0,0,0,0,  0,-3,0,0, 0,0,0,0, 0,0,0,0, 0,0,0,0];\\

From the Magma it can be caculated that the kissing number is indeed $146880$. We also know from Magma that the above two extremal even unimodualr lattice of dimension $32$ are isometric and they are actually the Barnes-Wall lattice ${\bf BW}_{32}$.\\

{\bf Remark 13.1.}  It seems possible to use more couples of binary linear $[16,k_1,4_2]$ and $[16, k_2,4]_2$ codes statisfying  $k_1+k_2=16$ such that the differences of the norms are at least "four" $\frac{9}{2}$'s or  two $\frac{9}{2}$'s plus one $9$, then we get  extremal even unimodular lattices of dimension $32$. We speculate that this would lead to some "new" extremal even unimodular lattices of dimension $32$ which are not in the list of fifteen Koch-Venkov extremal even unimodular lattices of dimension $32$.\\

{\bf 2nd}\\

Let ${\bf C}_4$ be the  linear ternary $[16, 8, 6]_3$ code with the generator matrix $({\bf I}_8|{\bf H}_8)$(\cite{CPS}), where ${\bf H}_8$ is the following Hadamard matrix.\\

$$
\left(
\begin{array}{cccccccccccccccccc}
1&1&1&1& 1&1&1&1\\
1&-1&-1&-1& 1&-1&1&1\\
1&-1&-1&1& -1&1&1&-1\\
1&-1&1&-1& 1&1&-1&-1\\
1&1&-1&1& 1&-1&-1&-1\\
1&-1&1&1& -1&-1&-1&1\\
1&1&1&-1& -1&-1&1&-1\\
1&1&-1&-1& -1&1&-1&1\\
\end{array}
\right)
$$

In our ternary construction we use the following $32$ dimensional lattice ${\bf L} \subset {\bf Z}^{16} \oplus {\bf Z}^{16} \oplus {\bf Z}^{16}$ defined by $x_{1j}+x_{2j}+x_{3j}=0$ for $j=1,...,16$ and $x_{i1}+\cdots+x_{i,16} \equiv 0$, $mod$ $2$. This is a lattice with minimum norm $4$ and volume $4 \cdot 3^8$. The above  ternary self-dual $[16, 8, 6]_3$ code ${\bf C}_4$ is used in the teranry construction. \\

We get a $32$ dimensional lattice ${\bf T'}_{32}^3$ with minimum norm $18$ and volume $4 \cdot 3^{32}$. The lattice is the union of $3^8$ translates of the lattice $3{\bf L}$ with leading vectors from the $3^8$ codewords in the ternary code ${\bf C}_3$. Let ${\bf x}=(\frac{3}{2}, ..., \frac{3}{2}) \in {\bf R}^{16}$ (sixteen $\frac{3}{2}$'s). Set ${\bf X}_1=({\bf x}, -{\bf x}, {\bf 0}), {\bf X}_2=({\bf 0}, {\bf x}, {\bf x})$ and ${\bf X}_3={\bf X}_1+{\bf X}_2$. The lattice ${\bf T}_{32}^3$ is the union of the four translates ${\bf T'}_{32}^3, {\bf X}_1+{\bf T'}_{32}^3$, ${\bf X}_2+{\bf T'}_{32}^3$ and ${\bf X}_3+{\bf T'}_{32}^3$.\\

{\bf Theorem 13.5.} {\em The minimum norm of the lattice ${\bf T}_{32}^3$ is $36$. The lattice $\frac{1}{3}{\bf T}_{32}^3$ is an unimodular even lattice with minimum norm $4$.}\\

{\bf Proof.} The point is that for any codeword of Hamming weight less than or equal to $12$ in the self-dual ternary $[16, 8, 6]_3$ code, the norm of the difference is at least $\frac{3 \cdot 12+ 4\cdot 9}{2}=36$. For any codeword in the self-dual ternary $[16, 8, 6]_3$ code with Hammming weight $15$, from the condition that the sum of each row is an even number the norm of the difference is at least $\frac{3 \cdot 15 +9}{2}+9=36$.\\

The following is a base of this lattice. It is easy to see that $3^2 \cdot 7$ is a factor of the order of the automorphism group of this lattice, since the automorphism group of this ternary $[16, 8, 6]_3$ group is obviously the subgroup of this lattice and the permutation group of three rows is obviously a subgroup of the automorphism group of this lattice. Thus if this is in the list of Koch-Venkov list of fifteen $32$ dimensional extremal even unimodular lattices (\cite{Nelist}), it has to be one of the lattice ${\bf BW}_{32}$, ${\bf Lambda(F)}$ and ${\bf Lambda(U)}$ lattices. We checked the "IsIsometric" in Magma with these three lattices, the Magma had no answer after three hours. Thus we conclude that here this 2nd lattice is a new $32$ dimensional extremal even unimodular lattice.\\

v1 := V![3/2,3/2,3/2,3/2, 3/2,3/2,3/2,3/2,3/2,3/2,3/2,3/2, 3/2,3/2, 3/2, 3/2,
-3/2,-3/2,-3/2,-3/2, -3/2,-3/2,-3/2,-3/2,-3/2,-3/2,-3/2,-3/2, -3/2,-3/2,-3/2,-3/2,
0,0,0,0, 0,0,0,0, 0,0,0,0, 0,0,0,0];\\

v2 :=V![0,0,0,0, 0,0,0,0, 0,0,0,0, 0,0,0,6, 0,0,0,0, 0,0,0,0, 0,0,0,0, 0,0,0,-6, 0,0,0,0, 0,0,0,0, 0,0,0,0, 0,0,0,0];\\

v3 :=V![0,3,0,0, 0,0,0,0, 0,0,0,0, 0,0,0,3, 0,-3,0,0, 0,0,0,0, 0,0,0,0, 0,0,0,-3, 0,0,0,0, 0,0,0,0, 0,0,0,0, 0,0,0,0];\\

v4 :=V![0,0,3,0, 0,0,0,0, 0,0,0,0, 0,0,0,3, 0,0,-3,0, 0,0,0,0, 0,0,0,0, 0,0,0,-3, 0,0,0,0, 0,0,0,0, 0,0,0,0, 0,0, 0,0];\\

v5 :=V![0,0,0,3, 0,0,0,0, 0,0,0,0, 0,0,0,3, 0,0,0,-3, 0,0,0,0, 0,0,0,0, 0,0,0,-3, 0,0,0,0, 0,0,0,0, 0,0,0,0, 0,0,0,0];\\

v6 :=V![0,0,0,0, 3,0,0,0, 0,0,0,0, 0,0,0,3,  0,0,0,0, -3,0,0,0, 0,0,0,0, 0,0,0,-3, 0,0,0,0, 0,0,0,0, 0,0,0,0, 0,0,0,0];\\

v7 :=V![0,0,0,0, 0,3,0,0, 0,0,0,0, 0,0,0,3,  0,0,0,0, 0,-3,0,0, 0,0,0,0, 0,0,0,-3, 0,0,0,0, 0,0,0,0, 0,0,0,0, 0,0,0,0];\\

v8 :=V![0,0,0,0, 0,0,3,0, 0,0,0,0, 0,0,0,3, 0,0,0,0, 0,0,-3,0, 0,0,0,0, 0,0,0,-3, 0,0, 0,0,0,0, 0,0,0,0, 0,0,0,0, 0,0];\\

v9 :=V![0,0,0,0, 0,0,0,3, 0,0,0,0, 0,0,0,3, 0,0,0,0, 0,0,0,-3, 0,0,0,0, 0,0,0,-3, 0,0, 0,0,0,0, 0,0,0,0, 0,0,0,0, 0,0];\\

v10 :=V![0,0,0,0, 0,0,0,0, 3,0,0,0, 0,0,0,3, 0,0,0,0,  0,0,0,0, -3,0,0,0, 0,0,0,-3, 0,0,0,0, 0,0,0,0, 0,0,0,0, 0,0,0,0];\\

v11 :=V![0,0,0,0, 0,0,0,0, 0,3,0,0, 0,0,0,3, 0,0,0,0,  0,0,0,0, 0,-3,0,0, 0,0,0,-3, 0,0,0,0, 0,0,0,0, 0,0,0,0, 0,0,0,0];\\

v12 :=V![0,0,0,0, 0,0,0,0, 0,0,3,0, 0,0,0,3, 0,0,0,0, 0,0,0,0, 0,0,-3,0,  0,0,0,-3, 0,0,0,0, 0,0,0,0, 0,0,0,0, 0,0,0,0];\\

v13 :=V![0,0,0,0, 0,0,0,0, 0,0,0,3, 0,0,0,3, 0,0,0,0, 0,0,0,0, 0,0,0,-3, 0,0, 0,-3, 0,0, 0,0,0,0, 0,0,0,0, 0,0,0,0, 0,0];\\

v14 :=V![0,0,0,0, 0,0,0,0, 0,0,0,0, 3,0,0,3, 0,0,0,0, 0,0,0,0,  0,0,0,0, -3,0,0,-3, 0,0,0,0, 0,0,0,0, 0,0,0,0, 0,0,0,0];\\

v15 :=V![0,0,0,0, 0,0,0,0, 0,0,0,0, 0,3,0,3, 0,0,0,0, 0,0,0,0,  0,0,0,0, 0,-3,0,-3, 0,0,0,0, 0,0,0,0, 0,0,0,0, 0,0,0,0];\\

v16 :=V![0,0,0,0, 0,0,0,0, 0,0,0,0, 0,0,3,3, 0,0,0,0, 0,0,0,0, 0,0,0,0, 0,0,-3,-3,  0,0,0,0, 0,0,0,0, 0,0,0,0, 0,0,0,0];\\

v17 := V![0,0,0,0, 0,0,0,0, 0,0,0,0, 0,0,0,0, -3/2,-3/2,-3/2,-3/2, -3/2,-3/2,-3/2,-3/2,-3/2,-3/2,-3/2,-3/2, -3/2,-3/2,-3/2,-3/2,
3/2,3/2,3/2,3/2, 3/2,3/2, 3/2,3/2,3/2,3/2,3/2,3/2, 3/2,3/2,3/2,3/2];\\

v18 :=V![0,0,0,0, 0,0,0,0, 0,0,0,0, 0,0,0,0, 0,0,0,0, 0,0,0,0, -3,0,0,0, 0,0,0,-3, 0,0,0,0, 0,0,0,0, 3,0,0,0, 0,0,0,3];\\

v19 :=V![0,0,0,0, 0,0,0,0, 0,0,0,0, 0,0,0,0, 0,0,0,0, 0,0,0,0, 0,3,0,0, 0,0,0,-9, 0,0,0,0, 0,0,0,0, 0,-3,0,0, 0,0,0,9];\\
 
v20 :=V![0,0,0,0, 0,0,0,0, 0,0,0,0, 0,0,0,0, 0,0,0,0, 0,0,0,0, 0,0,3,0, 0,0,0,-9, 0,0,0,0, 0,0,0,0, 0,0,-3,0, 0,0,0,9];\\

v21 :=V![0,0,0,0, 0,0,0,0, 0,0,0,0, 0,0,0,0, 0,0,0,0, 0,0,0,0, 0,0,0,3, 0,0,0,-9, 0,0,0,0, 0,0,0,0, 0,0,0,-3, 0,0,0,9];\\

v22 :=V![0,0,0,0, 0,0,0,0, 0,0,0,0, 0,0,0,0,  0,0,0,0, 0,0,0,0, 0,0,0,0, 3,0,0,-9,  0,0,0,0, 0,0,0,0, 0,0,0,0, -3,0,0,9];\\

v23 :=V![0,0,0,0, 0,0,0,0, 0,0,0,0, 0,0,0,0, 0,0,0,0, 0,0,0,0, 0,0,0,0, 0,3,0,-9, 0,0,0,0, 0,0,0,0, 0,0,0,0, 0, -3,0,9];\\

v24 :=V![0,0,0,0, 0,0,0,0, 0,0,0,0, 0,0,0,0, 0,0,0,0, 0,0,0,0, 0,0,0,0, 0,0,3,-9, 0,0,0,0, 0,0,0,0, 0,0,0,0, 0, 0,-3,9];\\

v25 :=V![-2,0,0,0, 0,0,0,0, 1,1,1,1, 1,1,1,1, 1,0,0,0, 0,0,0,0, 1,-2,-2,-2, -2,-2,-2,-2, 1,0,0,0, 0,0,0,0, -2,1,1,1, 1,1,1,1];\\

v26 :=V![0,-2,0,0, 0,0,0,0, 1,-1,-1,-1, 1,-1,1,1, 0,1,0,0, 0,0,0,0, 1,2,2,2, -2,2,-2,-2, 0,1,0,0, 0,0,0,0, -2,-1,-1,-1, 1,-1,1,1];\\

v27 :=V![0,0,-2,0, 0,0,0,0, 1,-1,-1,1, -1,1,1,-1, 0,0,1,0, 0,0,0,0, 1,2,2,-2, 2,-2,-2,2, 0,0,1,0, 0,0,0,0, -2,-1,-1,1, -1,1,1,-1];\\

v28 :=V![0,0,0,-2, 0,0,0,0, 1,-1,1,-1, 1,1,-1,-1, 0,0,0,1, 0,0,0,0, 1,2,-2,2, -2,-2,2,2, 0,0,0,1, 0,0,0,0, -2,-1,1,-1, 1,1,-1,-1];\\

v29 :=V![0,0,0,0, -2,0,0,0, 1,1,-1,1, 1,-1,-1,-1, 0,0,0,0, 1,0,0,0, 1,-2,2,-2, -2,2,2,2, 0,0,0,0, 1,0,0,0, -2,1,-1,1, 1,-1,-1,-1];\\

v30 :=V![0,0,0,0, 0,-2,0,0, 1,-1,1,1, -1,-1,-1,1, 0,0,0,0, 0,1,0,0, 1,2,-2,-2, 2,2,2,-2, 0,0,0,0, 0,1,0,0, -2,-1,1,1, -1,-1,-1,1];\\

v31 :=V![0,0,0,0, 0,0,-2,0, 1,1,1,-1, -1,-1,1,-1, 0,0,0,0, 0,0,1,0, 1,-2,-2,2, 2,2,-2,2, 0,0,0,0, 0,0,1,0, -2,1,1,-1, -1,-1,1,-1];\\

v32 :=V![0,0,0,0, 0,0,0,-2, 1,1,-1,-1, -1,1,-1,1, 0,0,0,0, 0,0,0,1, 1,-2,2,2, 2,-2,2,-2, 0,0,0,0, 0,0,0,1, -2,1,-1,-1, -1,1,-1,1];\\

{\bf 3rd}\\

Actually if the vector ${\bf x}$ is replaced by $(\pm \frac{3}{2},...,\pm\frac{3}{2}) \in {\bf R}^{16}$.... an extremal even unimodular lattice of dimension $32$ can be consctructed. The number $3^2 \cdot 7$ is a factor of the order of the automorphism of this lattice. Set ${\bf x}=(-\frac{3}{2}, \frac{3}{2},...,\frac{3}{2}) \in {\bf R}^{16}$. We get a $32$ dimensional extremal even unimodular lattice $\frac{1}{3}{\bf T}_{32}^4$.\\

The following is a base of this ${\bf T}_{32}^4$ . We asked Magma if this 3rd lattice is isometric to the 2nd one, the Magma had no answer after three hours. We conclude that the 3rd lattice and the 2nd are not isometric.\\

v1 := V![-3/2,3/2,3/2,3/2, 3/2,3/2,3/2,3/2,3/2,3/2,3/2,3/2, 3/2, 3/2, 3/2, 3/2,
3/2,-3/2,-3/2,-3/2, -3/2,-3/2,-3/2,-3/2,-3/2,-3/2,-3/2,-3/2, -3/2,-3/2,-3/2,-3/2,
0,0,0,0, 0,0,0,0, 0,0,0,0, 0,0,0,0];\\

v2 :=V![0,0,0,0, 0,0,0,0, 0,0,0,0, 0,0,0,6, 0,0,0,0, 0,0,0,0, 0,0,0,0, 0,0,0,-6, 0,0,0,0, 0,0,0,0, 0,0,0,0, 0,0,0,0];\\

v3 :=V![0,3,0,0, 0,0,0,0, 0,0,0,0, 0,0,0,3, 0,-3,0,0, 0,0,0,0, 0,0,0,0, 0,0,0,-3, 0,0,0,0, 0,0,0,0, 0,0,0,0, 0,0,0,0];\\

v4 :=V![0,0,3,0, 0,0,0,0, 0,0,0,0, 0,0,0,3, 0,0,-3,0, 0,0,0,0, 0,0,0,0, 0,0,0,-3, 0,0,0,0, 0,0,0,0, 0,0,0,0, 0,0, 0,0];\\

v5 :=V![0,0,0,3, 0,0,0,0, 0,0,0,0, 0,0,0,3, 0,0,0,-3, 0,0,0,0, 0,0,0,0, 0,0,0,-3, 0,0,0,0, 0,0,0,0, 0,0,0,0, 0,0,0,0];\\

v6 :=V![0,0,0,0, 3,0,0,0, 0,0,0,0, 0,0,0,3,  0,0,0,0, -3,0,0,0, 0,0,0,0, 0,0,0,-3, 0,0,0,0, 0,0,0,0, 0,0,0,0, 0,0,0,0];\\

v7 :=V![0,0,0,0, 0,3,0,0, 0,0,0,0, 0,0,0,3,  0,0,0,0, 0,-3,0,0, 0,0,0,0, 0,0,0,-3, 0,0,0,0, 0,0,0,0, 0,0,0,0, 0,0,0,0];\\

v8 :=V![0,0,0,0, 0,0,3,0, 0,0,0,0, 0,0,0,3, 0,0,0,0, 0,0,-3,0, 0,0,0,0, 0,0,0,-3, 0,0, 0,0,0,0, 0,0,0,0, 0,0,0,0, 0,0];\\

v9 :=V![0,0,0,0, 0,0,0,3, 0,0,0,0, 0,0,0,3, 0,0,0,0, 0,0,0,-3, 0,0,0,0, 0,0,0,-3, 0,0, 0,0,0,0, 0,0,0,0, 0,0,0,0, 0,0];\\

v10 :=V![0,0,0,0, 0,0,0,0, 3,0,0,0, 0,0,0,3, 0,0,0,0,  0,0,0,0, -3,0,0,0, 0,0,0,-3, 0,0,0,0, 0,0,0,0, 0,0,0,0, 0,0,0,0];\\

v11 :=V![0,0,0,0, 0,0,0,0, 0,3,0,0, 0,0,0,3, 0,0,0,0,  0,0,0,0, 0,-3,0,0, 0,0,0,-3, 0,0,0,0, 0,0,0,0, 0,0,0,0, 0,0,0,0];\\

v12 :=V![0,0,0,0, 0,0,0,0, 0,0,3,0, 0,0,0,3, 0,0,0,0, 0,0,0,0, 0,0,-3,0,  0,0,0,-3, 0,0,0,0, 0,0,0,0, 0,0,0,0, 0,0,0,0];\\

v13 :=V![0,0,0,0, 0,0,0,0, 0,0,0,3, 0,0,0,3, 0,0,0,0, 0,0,0,0, 0,0,0,-3, 0,0, 0,-3, 0,0, 0,0,0,0, 0,0,0,0, 0,0,0,0, 0,0];\\

v14 :=V![0,0,0,0, 0,0,0,0, 0,0,0,0, 3,0,0,3, 0,0,0,0, 0,0,0,0,  0,0,0,0, -3,0,0,-3, 0,0,0,0, 0,0,0,0, 0,0,0,0, 0,0,0,0];\\

v15 :=V![0,0,0,0, 0,0,0,0, 0,0,0,0, 0,3,0,3, 0,0,0,0, 0,0,0,0,  0,0,0,0, 0,-3,0,-3, 0,0,0,0, 0,0,0,0, 0,0,0,0, 0,0,0,0];\\

v16 :=V![0,0,0,0, 0,0,0,0, 0,0,0,0, 0,0,3,3, 0,0,0,0, 0,0,0,0, 0,0,0,0, 0,0,-3,-3,  0,0,0,0, 0,0,0,0, 0,0,0,0, 0,0,0,0];\\

v17 := V![0,0,0,0, 0,0,0,0, 0,0,0,0, 0,0,0,0, 3/2,-3/2,-3/2,-3/2, -3/2,-3/2,-3/2,-3/2,-3/2,-3/2,-3/2,-3/2, -3/2,-3/2,-3/2,-3/2,
-3/2,3/2,3/2,3/2, 3/2,3/2, 3/2,3/2,3/2,3/2,3/2,3/2, 3/2,3/2,3/2,3/2];\\

v18 :=V![0,0,0,0, 0,0,0,0, 0,0,0,0, 0,0,0,0, 0,0,0,0, 0,0,0,0, -3,0,0,0, 0,0,0,-3, 0,0,0,0, 0,0,0,0, 3,0,0,0, 0,0,0,3]; \\

v19 :=V![0,0,0,0, 0,0,0,0, 0,0,0,0, 0,0,0,0, 0,0,0,0, 0,0,0,0, 0,3,0,0, 0,0,0,-3, 0,0,0,0, 0,0,0,0, 0,-3,0,0, 0,0,0,3];\\
 
v20 :=V![0,0,0,0, 0,0,0,0, 0,0,0,0, 0,0,0,0, 0,0,0,0, 0,0,0,0, 0,0,3,0, 0,0,0,-3, 0,0,0,0, 0,0,0,0, 0,0,-3,0, 0,0,0,3];\\

v21 :=V![0,0,0,0, 0,0,0,0, 0,0,0,0, 0,0,0,0, 0,0,0,0, 0,0,0,0, 0,0,0,3, 0,0,0,-3, 0,0,0,0, 0,0,0,0, 0,0,0,-3, 0,0,0,3]; \\

v22 :=V![0,0,0,0, 0,0,0,0, 0,0,0,0, 0,0,0,0,  0,0,0,0, 0,0,0,0, 0,0,0,0, 3,0,0,-3,  0,0,0,0, 0,0,0,0, 0,0,0,0, -3,0,0,3];\\

v23 :=V![0,0,0,0, 0,0,0,0, 0,0,0,0, 0,0,0,0, 0,0,0,0, 0,0,0,0, 0,0,0,0, 0,3,0,3, 0,0,0,0, 0,0,0,0, 0,0,0,0, 0, -3,0,-3];\\

v24 :=V![0,0,0,0, 0,0,0,0, 0,0,0,0, 0,0,0,0, 0,0,0,0, 0,0,0,0, 0,0,0,0, 0,0,3,-3, 0,0,0,0, 0,0,0,0, 0,0,0,0, 0, 0,-3,3];\\

v25 :=V![-2,0,0,0, 0,0,0,0, 1,1,1,1, 1,1,1,1, 1,0,0,0, 0,0,0,0, 1,-2,-2,-2, -2,-2,-2,-2, 1,0,0,0, 0,0,0,0, -2,1,1,1, 1,1,1,1];\\

v26 :=V![0,-2,0,0, 0,0,0,0, 1,-1,-1,-1, 1,-1,1,1, 0,1,0,0, 0,0,0,0, 1,2,2,2, -2,2,-2,-2, 0,1,0,0, 0,0,0,0, -2,-1,-1,-1, 1,-1,1,1];\\

v27 :=V![0,0,-2,0, 0,0,0,0, 1,-1,-1,1, -1,1,1,-1, 0,0,1,0, 0,0,0,0, 1,2,2,-2, 2,-2,-2,2, 0,0,1,0, 0,0,0,0, -2,-1,-1,1, -1,1,1,-1];\\

v28 :=V![0,0,0,-2, 0,0,0,0, 1,-1,1,-1, 1,1,-1,-1, 0,0,0,1, 0,0,0,0, 1,2,-2,2, -2,-2,2,2, 0,0,0,1, 0,0,0,0, -2,-1,1,-1, 1,1,-1,-1];\\

v29 :=V![0,0,0,0, -2,0,0,0, 1,1,-1,1, 1,-1,-1,-1, 0,0,0,0, 1,0,0,0, 1,-2,2,-2, -2,2,2,2, 0,0,0,0, 1,0,0,0, -2,1,-1,1, 1,-1,-1,-1];\\

v30 :=V![0,0,0,0, 0,-2,0,0, 1,-1,1,1, -1,-1,-1,1, 0,0,0,0, 0,1,0,0, 1,2,-2,-2, 2,2,2,-2, 0,0,0,0, 0,1,0,0, -2,-1,1,1, -1,-1,-1,1];\\

v31 :=V![0,0,0,0, 0,0,-2,0, 1,1,1,-1, -1,-1,1,-1, 0,0,0,0, 0,0,1,0, 1,-2,-2,2, 2,2,-2,2, 0,0,0,0, 0,0,1,0, -2,1,1,-1, -1,-1,1,-1];\\

v32 :=V![0,0,0,0, 0,0,0,-2, 1,1,-1,-1, -1,1,-1,1, 0,0,0,0, 0,0,0,1, 1,-2,2,2, 2,-2,2,-2, 0,0,0,0, 0,0,0,1, -2,1,-1,-1, -1,1,-1,1];\\

As the 2nd lattice $3^2 \cdot 7$ is a factor of the order of the automorphism group of this 3rd lattice.Thus if this is in the list of Koch-Venkov list of fifteen $32$ dimensional extremal even unimodular lattices (\cite{Nelist}), it has to be one of the lattice ${\bf BW}_{32}$, ${\bf Lambda(F)}$ and ${\bf Lambda(U)}$ lattices. We checked the "IsIsometric" in Magma with these three lattices, the Magma had no answer after three hours. Thus we conclude here this 3rd lattice is another new $32$ dimensional extremal even unimodular lattice.\\

{\bf 4th}

In this $32$ dimensional extremal even unimodular lattice, we use the same ${\bf T'}_{32}^3$ as in the third one, the vectors ${\bf X}_1, {\bf X}_2$ and ${\bf X}_3$ are replaced by ${\bf X'}_1=({\bf x'}, -{\bf x'},{\bf 0}), {\bf X'}_2=({\bf 0}, {\bf x'}, -{\bf x'})$ and ${\bf X'}_3={\bf X'}_1+{\bf X'}_2$, where ${\bf x'}=(\frac{3}{2}, ..., \frac{3}{2}, 0,...,0) \in {\bf R}^{16}$ (eight $\frac{3}{2}$'s at the first eight positions and eight $0$'s at the last eight positions). The $32$ dimensional lattice ${\bf T}_{32}^4$ is the union of the four translates ${\bf T'}_{32}^3$, ${\bf X'}_1+{\bf T'}_{32}^3$, ${\bf X'}_2+{\bf T'}_{32}^3$ and ${\bf X'}_3+{\bf T'}_{32}^3$.\\

{\bf Theorem 13.6} {\em The minimum norm of the lattice ${\bf T}_{32}^4$ is $36$ and the lattice $\frac{1}{3} {\bf T}_{32}^4$ is a $32$ dimensional extremal even lattice of minimum norm $4$.}\\

{\bf Proof.} First of all it is obvious ${\bf H}_8 {\bf H}_8=-{\bf I}_8$ and $(-{\bf H}_8, {\bf I}_8)$ is also a generator matrix of the above $[16, 8, 6]_3$ ternary self-dual code. Then for a weight six codeword in this code, there are at  most four non-zero coordinate positions in the set $\{1,2,...,8\}$. The norms of the differences of ${\bf X'}_i$'s  with the lattice vectors in ${\bf T'}_{32}^3$ are at least $\frac{3 \cdot 4+ 9 \cdot 4}{2}+12=36$. For  codewords with Hamming weights $12,15$ in this ternary code, the norm of the differences are at least $\frac{3 \cdot 8}{2}+24=36$. In the Hamming weight $9$ case, we note if the norm of the difference at the position set $\{1,2,...,8\}$ is $\frac{3 \cdot 8}{2}=12$, the corresponding lattice vectors in this position set $\{1,2,...,8\}$ satisfying $\Sigma_{i=1}^8 t_i \equiv 0$ $mod$ $2$, where $t_i$'s are the coordinates of the lattice vector in ${\bf T'}_{32}^3$, thus the only non-zero column of this lattice vector in the position set $\{1,2,...,8\}$ has to be the form $\pm(2,-4,2)^{\tau}$ (or permutations). The norm is at least $12+24=36$. If the norm of this difference at the position set is not $12$,  thare are at least two additions of $9$ in the norm. Thus the norm of the differences are at least $36$.\\

We cannot decide now if this 4th $32$ dimensional extremal even unimodular lattice is new.\\

\section{A new 48 dimensional extremal even unimodular  lattice}

In this section we indicate how our method leads to $48$ dimensional extremal unimodular even lattices with minimum norm $6$ naturally. This gives a possibility that a few {\bf new} such lattices can be constructed. For the recent progress we refer to \cite{Nebe1}.\\

In our ternary construction we use the following $48$ dimensional lattice ${\bf L} \subset {\bf Z}^{24} \oplus {\bf Z}^{24} \oplus {\bf Z}^{24}$ satisfying $x_{1j}+x_{2j}+x_{3j}=0$ for $j=1,...,24$ and $x_{i1}+\cdots+x_{i24} \equiv 0$ $mod$ $2$ for $i=1,2,3$.  A ternary self-dual $[24, 12, 6]_3$ or $[24, 12, 9]_3$ code ${\bf C}$ is used in the teranry construction we get a $48$ dimensional lattice ${\bf T'}_{48}$ with minimum norm $36$ and volume $ 4\cdot 3^{48}$. The lattice is the union of $3^{12}$ translates of the lattice $3{\bf L}$ with leading vectors from the $3^{12}$ codewords in the ternary code ${\bf C}$.\\

We need to use $2^{24}$ translates of this lattice ${\bf T'}_{48}$ according to the binary Golay $[24, 12, 8]_2$ code. The union of these $2^{24}$ translates is a new $48$ dimensional lattice ${\bf T''}_{48}$ with the volume $\frac{4\cdot 3^{48}}{2^{24}}=\frac{3^{48}}{2^{22}}$ with the minimum norm $27$.  Without loss of generality suppose there are two weight $12$ codewords in the binary Golay $[24, 12, 8]_2$ code (a "{\em duum"} as named in page 279 of \cite{CS1}) supported at the position sets $\{1,...,12\}$ and $\{13,...,24\}$. We use a base for the binary Golay $[24, 12, 8]_2$  code with the first two element these two weight $12$ codewords. For  the first weight $12$ codewords in the base, set  ${\bf x}_c$ as $(-\frac{3}{2}, -\frac{3}{2}, -\frac{3}{2},\frac{3}{2},...,\frac{3}{2})$, which has three $-\frac{3}{2}$'s and nine $\frac{3}{2}$''s. For the second twelve positions, ${\bf x}_c$ as $(-\frac{3}{2}, \frac{3}{2}, \frac{3}{2},\frac{3}{2},...,\frac{3}{2})$, which has one $-\frac{3}{2}$ and eleven  $\frac{3}{2}$''s. For any other element $c$  in the base whose support positions are $i_1,...,i_t$, set ${\bf x}_c$ the vector in ${\bf R}^{24}$ with $\frac{3}{2}$ at these positions $i_1,..,i_t$. Let ${\bf X}_{c, 1}=({\bf x}_c, -{\bf x}_c,{\bf 0})$ and ${\bf X}_{c, 2}=({\bf 0}, {\bf x}_c, -{\bf x}_c)$. For any codeword $c$ in the binary Golay code which is expressed by a linear combination of codewords in the base, we just use the same linear conbination to get ${\bf X}_{c,1}$ and ${\bf X}_{c,2}$. Then we have $2^{12+12}$ such vectors.\\

The main point here is that the minimum norm of the differences in the translates is at least $27$.\\

We recall that in a ternary self-dual code, the Hamming weights of all codewords can be divided by $3$, that implies there are only weights $6, 9, 12, 15$, $18, 21$ and $24$ codewords in ${\bf C}$.\\

{\bf Theorem 12.1.} {\em If the following conditions are satisfied the minimum norm of the lattice ${\bf T''}_{48}$ is $27$.\\
1) For any weight $6$ codeword in the ternary code ${\bf C}$, its support is not in the support of a weight $8$ codeword in the binary Golay code;\\
2) For any weight $9$ codeword in the ternary code ${\bf C}$, its support does not contain the support of a weight $8$ codeword in the binary Golay code;\\ 
3) For any weight $12$ codeword in the ternary code ${\bf C}$, its support does not equal to the support of a weight $12$ codeword in the binary Golay code and does not equal to the union of the supports of two weight $8$ codewords in the binary Golay code}\\

{\bf Proof.} For any ${\bf X}_{c,1}-{\bf X}_{c',2}-{\bf e}$, where  $e$ is a lattice vector in ${\bf T'}_{48}$, if ${\bf e}$ is in $3{\bf L}$, the norm of this difference is at least $36$. If ${\bf e}$ is not in $3{\bf L}$ and it is from a cordword $c''$ of weights
 $6, 9, 12, 15,18, 21, 24$ in the ternary code ${\bf C}$, the norm of the difference is at least $6t_1+\frac{3t_2}{2}+\frac{9t_3}{2}$, where $t_1$ is the number of the positions of the support of $c''$ outside the union of the supports of $c$ and $c'$, $t_2$ is the number of positions of the support of $c''$ in the union of the supports of $c$ and $c'$, and $t_3$ is the number of the positions of the union of supports of $c$ and $c'$ outside the support of $c''$. It is clear that the union of  the supports of two codewords in the binary Golay code has $8$ or $12$ or $14$ or $16$ or $18$ or $20$ or $22$ or $24$ elements.\\

1) If the weight of $c''$ is at least $15$, since the union of the support $c$ and $c'$ cannot have odd number of positions (the binary Golay code is self-orthogonal), thus the norm of the difference is at least $\frac{3 \cdot 15 +9}{2}=\frac{3 \cdot 14}{2}+6=27$.\\

2) If the weight of $c''$ is $12$ and not supported at one of the above {\em "duum"}, the norm of the difference is at least $\frac{3 \cdot 12}{2}=18$, if this possibility does not happen, the increasing of norm of the difference is at least $9$. Thus from the condition 3) the norm of the difference is at least $18+9=27$. For the two weight $12$ codewords supported at  the duum, since there is only one $-\frac{3}{2}$ in ${\bf x}_c$, the minimum norm of the difference is at least $36$ from the same argument as in section 6. \\

3) If the weight of $c''$ is $9$, the norm of the difference is at least $\frac{3 \cdot 8}{2}+6=18$, if this possibility cannot happen, the increasing of norm of the difference is at least $9$. Thus from the condition 2) the norm of the difference is at least $18+9=27$.\\

4) If the weight of $c''$ is $6$, the norm of the difference attains the minimum on when   $c=c'$ is  a weight $8$ codeword in the binary Golay code. The minimum norm of the difference is $\frac{3 \cdot 6+2 \cdot 9}{2}=18$, if this possibility does not  happen, the increasing of norm of the difference is at least $9$. Thus from the condition 1) the norm of the difference is at least $18+9=27$.\\

On the other hand the norm of the difference of the vector $(\frac{3}{4}, -\frac{3}{4}, 0)^{\tau}$ with $(1,-2,1)^{\tau}$ is bigger than $\frac{9}{16}+\frac{9}{16}=\frac{9}{8}$. Thus the norm of the difference involving ${\bf Y}_1, {\bf Y}_2$ is at least $24 \cdot \frac{9}{8}=27$. The conclusion is proved.\\

From the classification of the paper \cite{Leon}, there are many inequivalent $[24, 12, 6]_3$ ternary code and exactly two inequivalent $[24, 12, 9]_3$ ternary code, if the positions can be arranged so that the conditions 1) 2) and 3) in Theorem 12.1 are satisfied, we do get an extremal even unimodular $48$ dimensional lattice with the minimum norm $6$. It would be difficult to prove these possible lattice are distinct.\\

The following result about the supports of codewords in the bianry Golay $[24, 12, 8]_2$ code is useful in our construction.\\

{\bf Proposition 12.2.} {\em Suppose there are two weight $12$ codewords supported at the position sets $\{1,...,12\}$ and $\{13,...,24\}$ in the binary Golay $[24, 12, 8]_2$ code. \\
1) Then any weight $8$ codeword in the binary Golay code is supported at two (or six) positions in the set $\{1,...,12\}$ and six (or two) positions in the set $\{13,...,24\}$ ($(2,6)$ or $(6,2)$ type), or is supported at four positions at both of these sets.\\
2) Any weight $12$ codeword in the binary Golay code is supported at four (or eight) positions in the set $\{1,...,12\}$ and eight (or four) positions in the set $\{13,...,24\}$ ($(4,8)$ or $(8,4)$ type), or is supported at six positions at both of these sets.}\\is 

{\bf Proof.} From Theorem 15 in page 280 of \cite{CS1}, the sub-group of the automorphism group preserving these two weight $12$ codewords is the Mathieu group ${\bf M}_{12}$. It  acts on the codewords of the binary Golay code naturally and it is $5$ transitive. Thus we have $2 \cdot 2 \cdot \displaystyle{12 \choose 2}=264$ type $(2,6)$ or $(6,2)$ weight $8$ codewords. On the other hand since ${\bf M}_{12}$ acts $5$ transitively, there are at least $\displaystyle{12 \choose 4}=495$ type $(4,4)$ weight $8$ codewords. These are totally $264+495=759$ weight $8$ codewords in the binary Golay code. We note that all type $(4, 4)$ weight $8$ codewords forms one orbit under the action of the group of the Mathieu group ${\bf M}_{12}$.\\

From the difference with the type $(4, 4)$ weight $8$ codewords, we have at least $495+495=990$ type $(4,8)$ or $(8,4)$ weight $12$ codewords. On the other hand the intersetion set of two weight $12$ codewords in the binary Golay code has $0$ or $4$ or $6$ or $8$ positions. Thus we have $2576-990-2=1584$ type $(6, 6)$ weight $12$ codewords in the binary Golay code. These $1584$ weight $12$ codewords forms $12$ orbits under the action of ${\bf M}_{12}$.\\

We will take the ternary $[24, 12, 6]_3$ code ${\bf C}$ as the direct sum of two copies of the extended ternary $[12, 6,, 6]_3$ Golay code which are supported in the sets $\{1,...,12\}$ and $\{13,...,24\}$. From Proposition 12.2, it is clear the condition 2) about the weight $9$ codewords in ${\bf C}$ is satisfied. \\

Since $\displaystyle{12 \choose 6}=132 \cdot 7$, all type $(6, 6)$ sets with $12$ elements (the subsets of $\{1,...,24\}$ with $12$ elements whose intersection sets with $\{1,...,12\}$ and $\{13,...,24\}$ have $6$ elements) forms $7 \cdot 7 \cdot 132 =49 \cdot 132=6468$ orbits under the action of the Mathieu group ${\bf M}_{12}$. Now we count the number of the orbits of type $(6, 6)$ subsets in the set $\{1,...,24\}$ which is of the following form (the subsets which attain the minimum norm of difference $18$ in the proof of Theorem 12.1).\\
1) One of the six-element-subset in $\{1,...,12\}$ or $\{13,...,24\}$ is in the support of a type $(2,6)$ or $(6,2)$ weight $8$ codeword in the binary Golay code;\\
2) The supports of type $(6,6)$ weight $12$ codewords in the binary Golay code;\\
3) The union of the supports of two weight $8$ codewords in the binary Golay code.\\

In 1) we have $1 \cdot 7+7 \cdot 1-1=13$ orbits, in 2) we have $12$ orbits. In the case 3), the difference of two weight $8$ codewords in the binary code is a weight $8$ codeword. The intersection of the supports of these two weight $8$ codewords is a set with $4$ elements. From Proposition 12.2, both weight $8$ codewords are\\
3a) type $(4,4)$ or\\
3b)one is type $(2,6)$ and another is type $(6,2)$.\\

It is clear that there only one orbit in case 3b). In case 3a), since the complementary set of $\{1,...,12\}$ of a four-element-set has eight elements. There are at most $\displaystyle{8 \choose 2}=28$ orbits. Totally we have $13+12+1+28=54$ orbits. Thus if we take the supports of one weight $6$ codeword in the first copy of the ternay $[12, 6, 6]_3$ code and another weight $6$ codeword in the second copy of ternary $[12, 6, 6]_3$ code not these $54$ orbits, the conditions of Theorem 12,2 are satisfied. We get a $48$ dimensional lattice  ${\bf T''}_{48}$ with the minimum norm $27$ and the volume $\frac{3^{48}}{2^{22}}$. \\

 Let ${\bf y}=(\frac{3}{4}, -\frac{3}{4},...,-\frac{3}{4},...,-\frac{3}{4},..., \frac{3}{4},...,-\frac{3}{4},...,\frac{3}{4}) \in {\bf R}^{24}$, which has nine $\frac{3}{4}$'s and three $-\frac{3}{4}$ at the first twelve positons, eleven $\frac{3}{4}$'s and one $-\frac{3}{4}$ at the last twelve positins. Set ${\bf Y}_1=({\bf y}, -{\bf y}, {\bf 0})$ and ${\bf Y}_2=({\bf 0}, {\bf y}, -{\bf y})$. The lattice ${\bf T}_{48}$ is the union of $2^{2}$ translates of ${\bf T''}_{48}$ with these $2^{2}$ leading vectors.\\

{\bf Theorem 12.3.} {\em The minimum norm of this $48$ dimensional lattice is $27$. The lattice $\frac{\sqrt{2}}{3}{\bf T}_{48}$ is an integral unimodular even lattice with the minimum norm $6$. }\\

{\bf Proof.} From the formation of the vectors ${\bf Y}_i$, $i=1,2$, it is clear that there are four $(-\frac{3}{4}, \frac{9}{4}, -\frac{3}{2})^{\tau}$'s if we want the part corresponding to the $\frac{3}{4}$'s to be the form $(\frac{3}{4},\frac{3}{4}, -\frac{3}{2})^{\tau}$. Then the norms of the difference for the columns $(\frac{3}{4},\frac{3}{4}, -\frac{3}{2})$ are at least $(\frac{-1}{4})^2+(\frac{-1}{4})^2+(\frac{1}{2})^2=\frac{3}{8}$. The nomrs of the differences for these four columns  $(-\frac{3}{4}, \frac{9}{4}, -\frac{3}{2})^{\tau}$'s are at least $(\frac{5}{4})^2+(\frac{-7}{4})^2+(\frac{1}{2})^2=\frac{39}{8}=\frac{3}{8}+\frac{9}{2}$. Thus totally the minimum norm is at least $24 \cdot \frac{3}{8}+4\cdot \frac{9}{2}=9+18=27$.\\

We denote the support vector of $3+1$ positions of $-\frac{3}{4}$'s as ${\bf x} \in {\bf F}_2^{24}$, if $g_i \in {\bf M}_{12}, i=1,2$ satisfying that ${\bf x}-g_i{\bf x}$ is in the binary linear Golay code, it is clear ${\bf x}-g_1g_2{\bf x}$ is in the binary linear Golay code. In fact $g_1^{-1}{\bf x}-{\bf x}+{\bf x}-g_2{\bf x}=g_1^{-1}{\bf x}-g_2{\bf x}$ is in the binary linear Golay code. Thus the set of elements in the Mathieu group ${\bf M}_{12}$ such that ${\bf x}-g{\bf x}$ is in the binary Golay $[24,12, 8]_2$ code is a subgroup ${\bf G}$. The set of elements in ${\bf M}_{12}$ fixing the set of three positions of $-\frac{3}{4}$ in $\{1,2,...,12\}$ and  fixing the one position of the $-\frac{3}{4}$ in $\{13,14,...,24\}$ is a subgroup of ${\bf G}$. From section 1.5 of Chapter 10 and Theorem 15 in page 280 of \cite{CS1}, it is easy to find an order five element in this subgroup.\\

It is clear that ${\bf G}$ is a sub-group of the automorphism group of the lattice ${\bf T}_{48}$.\\

{\bf Theorem 12.4.}{\em The order of the subgroup ${\bf G}$ of the  Mathieu group ${\bf M}_{12}$ can be divided by $55$. This ${\bf G}$ is a subgroup  of the automorphism group of the $48$ dimensional extremal even unimodular lattice $\frac{\sqrt{2}}{3} {\bf T}_{48}$.}\\

{\bf Proof.} By a suitable arrangement of the eight positions of an type $(6,2)$ octad in the set $\{1,2,...,24\}$  we can make that any prescribed element in ${\bf M}_{12}$ is an element of ${\bf G}$. Then ${\bf G}$ can have an element of ${\bf M}_{12}$ with the order $11$. Thus the order of ${\bf G}$ can be divieded by $55$.  From the construction we can verify that $\frac{\sqrt{2}}{3} {\bf T}_{48}$ is an integral even unimodular lattice.\\

From Theorem 12.4 this $48$ dimensional extremal even unimodular  lattice is not the lattice ${\bf P}_{48p}$, ${\bf P}_{48q}$ and ${\bf P}_{48n}$ (\cite{CS1,Nebe1}. This  lattice is not the new Nebe's $48$ dimensionla extremal even unimodular  lattice ${\bf P}_{48m}$ in \cite{Nebe1,Nebe3} since the order of $Aut({\bf P}_{48m})$ is $1200$ (\cite{Nebe1}).\\

{\bf Remark 12.1.} 1) It is interesting to compare Theorem 12.4 with Theorem 5.8 of \cite{Nebe2}.\\

2) We speculate that $C$ could be many $[24, 12, 6]_3$ ternary codes in the paper \cite{Leon} and more new $48$ extremal even unimordular lattices could be constructed from Theorem 12.1.\\

{\bf Acknowledgement.} The author is grateful to Henry Cohn for his help and encouragement. This work was supported by National Natural Science Foundation of China Grants 11061130539 and  11371138.\\

\end{document}